\documentclass[12pt,leqno]{amsart}

\usepackage{amssymb, amsmath, amsfonts, latexsym}

\newtheorem{thm}[equation]{Theorem}
\newtheorem{lemma}[equation]{Lemma}
\newtheorem{prop}[equation]{Proposition}
\newtheorem{cor}[equation]{Corollary}

\theoremstyle{definition}
\newtheorem{definition}[equation]{Definition}

\newtheorem{assumptions}[equation]{Assumptions}

\theoremstyle{remark}
\newtheorem{remark}[equation]{Remark}

\newtheorem{examples}[equation]{Examples}

\numberwithin{equation}{section}

\begin{document}
%
%
\renewcommand{\P}{\mathbb{P}}
\newcommand{\E}{\mathbb{E}}
\newcommand{\1}{\textbf{1}}

\newcommand{\eqdef}{\stackrel{\vartriangle}{=}}
\newcommand{\eqlaw}{\stackrel{Law}{=}}
\newcommand{\dom}{\mathrm{dom\,}}
\newcommand{\icordom}{\mathrm{icordom\,}}
\newcommand{\cl}{\mathrm{cl\,}}
\newcommand{\cv}{\mathrm{cv}}
\newcommand{\inter}{\mathrm{int\,}}
\newcommand{\ls}{\mathrm{ls}\,}

\newcommand{\push}{\diamond}

\newcommand{\lsc}{lower semicontinuous}

\newcommand{\ldp}{LDP}
\newcommand{\ld}{LD}
\newcommand{\kn}{{(k,n)}}
\newcommand{\knldp}{$\kn$-LDP}
\newcommand{\knld}{$\kn$-LD}
\newcommand{\nldp}{$n$-LDP}
\newcommand{\kldp}{$k$-LDP}

\newcommand{\limn}{\lim_{n\rightarrow\infty}}
\newcommand{\Glimn}{\Gamma\textrm{-}\lim_{n\rightarrow\infty}}
\newcommand{\lin}{\liminf_{n\rightarrow\infty}}
\newcommand{\lsn}{\limsup_{n\rightarrow\infty}}
\newcommand{\nlog}{\frac{1}{n}\log \mathbb{P}}

\newcommand{\limk}{\lim_{k\rightarrow\infty}}
\newcommand{\Glimk}{\Gamma\textrm{-}\lim_{k\rightarrow\infty}}
\newcommand{\lik}{\liminf_{k\rightarrow\infty}}
\newcommand{\lsk}{\limsup_{k\rightarrow\infty}}
\newcommand{\klog}{\frac{1}{k}\log P^k_z}

\newcommand{\limkn}{\limk\limn}
\newcommand{\likn}{\lik\lin}
\newcommand{\lskn}{\lsk\lsn}
\newcommand{\knlog}{\frac{1}{kn}\log \P}

\newcommand{\lima}{\lim_{\alpha\rightarrow\infty}}
\newcommand{\Glima}{\Gamma\textrm{-}\lim_{\alpha\rightarrow\infty}}

\newcommand{\Rd}{\mathbb{R}^d}
\newcommand{\Rdd}{\mathbb{R}^{2d}}
\newcommand{\X}{\mathcal{X}}
\newcommand{\PX}{\mathcal{P}_\X}
\newcommand{\Z}{\mathcal{Z}}
\newcommand{\PZ}{\mathcal{P}_\Z}
\newcommand{\ZX}{\Z\!\times\!\X}
\newcommand{\PZX}{\mathcal{P}_{\Z\X}}
\newcommand{\CX}{C_\X}
\newcommand{\UX}{U_{\X}}
\newcommand{\CZ}{C_\Z}
\newcommand{\CZX}{C_{\Z\X}}
\newcommand{\CR}{C_{\Rd}}
\newcommand{\CRR}{C_{\mathbb{R}^{2d}}}
\newcommand{\UZX}{U_{\Z\X}}
\newcommand{\PR}{\mathcal{P}_{\Rd}}
\newcommand{\PRR}{\mathcal{P}_{\mathbb{R}^{2d}}}
\newcommand{\AC}{\mathcal{AC}}

\newcommand{\cY}{{c^Y}}

\newcommand{\zni}{{z_{n,i}}}
\newcommand{\Xni}{{X_{n,i}}}
\newcommand{\Xkni}{{X^k_{n,i}}}
\newcommand{\Ykni}{{X^k_{n,i}}(1)}
\newcommand{\Wkni}{{W^k_{n,i}}}

\newcommand{\Lkn}{L^k_n}
\newcommand{\Mkn}{M^k_n}
\newcommand{\Nkn}{N^k_n}
\newcommand{\Kkn}{K^k_n}
\newcommand{\LWkn}{L^{W,k}_n}
\newcommand{\NWkn}{N^{W,k}_n}

\newcommand{\pk}{\pi^k}

\newcommand{\un}{\frac{1}{n}}
\newcommand{\IX}{\int_\X}
\newcommand{\IZ}{\int_\Z}
\newcommand{\IZX}{\int_{\ZX}}

\newcommand{\supx}{\sup_{x\in\X}}
\newcommand{\infx}{\inf_{x\in\X}}

\newcommand{\LDDa}{LDDa}
\newcommand{\LDDb}{LDDb}

\newcommand{\infX}{\inf_{x\in X}}
\newcommand{\infY}{\inf_{y\in Y}}
\newcommand{\infV}{\inf_{x\in V}}
\newcommand{\supY}{\sup_{y\in Y}}

\newcommand{\xy}{\langle x,y\rangle}



 \title{A large deviation approach to optimal transport}

 \author{Christian L\'eonard}
 \address{(Christian L\'eonard) Modal-X, Universit\'e Paris 10. B\^at. G, 200 av. de la R\'epublique. 92001 Nanterre Cedex, France}
 \address{(Christian L\'eonard) CMAP, \'Ecole Polytechnique. 91128 Palaiseau Cedex, France}
 \email{christian.leonard@polytechnique.fr}

\date{December 2005}



\begin{abstract}
A probabilistic method for solving the Monge-Kantorovich mass
transport problem on $\mathbb{R}^d$ is introduced. A system of
empirical measures of independent particles is built in such a way
that it obeys a doubly indexed large deviation principle with an
optimal transport cost as its rate function. As a consequence, new
approximation results for the optimal cost function and the
optimal transport plans are derived. They follow from the
$\Gamma$-convergence of a sequence of normalized relative
entropies toward the optimal transport cost. A wide class of cost
functions including the standard power cost functions $|x-y|^p$
enter this framework.
\end{abstract}

\maketitle


\section{Introduction}

This paper introduces a probabilistic method for solving the
Monge-Kantorovich mass transport problem.
\subsection{The Monge-Kantorovich problem}
Let $\mu$ and $\nu$ be two probability measures on $\Rd$ seen as
mass distributions. One wants to transfer $\mu$ to $\nu$ with a
minimal cost, given that transporting a unit mass from $x_0$ to
$x_1$ costs $c(x_0,x_1).$ This means that one searches for a
transport plan $x_1=T(x_0)$ such that the image measure
$T\push\mu$ is $\nu$ and $\int_{\Rd}c(x_0,T(x_0))\,\mu(dx_0)$ is
minimal. This problem was addressed by G.~Monge \cite{Monge} at
the eighteenth century. In the 40's, L.V.~Kantorovich
\cite{Kanto42}, \cite{Kanto48} proposed a relaxed version of Monge
problem by allowing each cell of mass at $x_0$ to crumble into
powder so that it can be tranfered to several $x_1$'s. In
mathematical terms, one searches for a probability measure $\rho$
on $\mathbb{R}^d\!\times\!\mathbb{R}^d$ whose marginal measures
$\rho_0(dx_0)=\rho(dx_0\times \Rd)$ and
$\rho_1(dx_1)=\rho(\Rd\times dx_1)$ satisfy $\rho_0=\mu$ and
$\rho_1=\nu$ and such that
$\int_{\mathbb{R}^d\!\times\!\mathbb{R}^d}c(x_0,x_1)\,\rho(dx_0dx_1)$
is minimal. Let us denote $\PR$ and $\PRR$ the sets of all
probability measures on $\Rd$ and
$\mathbb{R}^d\!\times\!\mathbb{R}^d.$ For each $\mu$ and $\nu$ in
$\PR,$ we face the optimization problem
\begin{equation}\label{eq-MK}
    \textrm{minimize\quad}
    \int_{\mathbb{R}^d\!\times\!\mathbb{R}^d}c(x_0,x_1)\,\rho(dx_0dx_1)
    \textrm{\quad subject to\quad}\rho\in\Pi(\mu,\nu) \tag{MK}
\end{equation}
where the cost function
$c:\mathbb{R}^d\!\times\!\mathbb{R}^d\mapsto [0,+\infty]$ is
assumed to be measurable and
\begin{equation*}
    \Pi(\mu,\nu)=\{\rho\in\PRR; \rho_0=\mu, \rho_1=\nu\}
\end{equation*}
is the set of all probability measures on
$\mathbb{R}^d\!\times\!\mathbb{R}^d$ with marginals $\mu$ and
$\nu.$ This problem is called the Monge-Kantorovich mass transport
problem. Monge problem corresponds to the transport plans
$\rho(dx_0dx_1)=\mu(dx_0)\delta_{T(x_0)}(dx_1)$ where $\delta$
stands for the Dirac measure. Kantorovich's relaxation procedure
embeds Monge's nonlinear problem in the linear programing problem
(MK).
\\
The value of (MK) is the \emph{transportation cost} defined for
all $\mu$ and $\nu$ in $\PR$ by
\begin{equation}\label{eq-12}
\mathcal{T}_c(\mu,\nu)
  :=\inf_{\rho\in
  \Pi(\mu,\nu)}\int_{\mathbb{R}^d\!\times\!\mathbb{R}^d}c(x_0,x_1)\,\rho(dx_0dx_1).
\end{equation}
The special cost function $c_p(x_0,x_1)=|x_1-x_0|^p$ with $p\geq
1,$ leads to the Wassertein metric
$\mathcal{T}_{c_p}^{1/p}(\mu,\nu).$

\subsection{Which large deviations?}
As the title of the paper indicates, our probabilistic approach of
Monge-Kantorovich problem is in terms of large deviations. One can
interpret $\mu$ and $\nu$ respectively as the distributions of the
initial and final random positions $X_0$ and $X_1$ of a random
process $(X_t)_{0\leq t\leq 1}.$ In the present paper, only the
couple of initial and final positions $(X_0,X_1)$ is considered.

Our aim is to obtain a Large Deviation Principle (LDP) in $\PR$
the rate function of which is $\nu\mapsto \mathcal{T}_c(\mu,\nu)$
where $\mu$ is fixed. The definition of a LDP is recalled at
(\ref{eq-ldp}). General cost functions will be considered in the
article but for the sake of clarity, in this introductory section
our procedure is described in the special case of the quadratic
cost funtion $c(x_0,x_1)=|x_1-x_0|^2/2.$  For each integer $k\geq
1,$ take a system of $n$ independent random couples
$(\Xkni(0),\Xkni(1))_{1\leq i\leq n}$ which is described as
follows. For each $i,$ the initial position $\Xkni(0)=\zni$ is
deterministic and the final position is
$$\Xkni(1)=\zni +Y_i/\sqrt{k} $$ where the $Y_i$'s are independent
copies of a standard normal vector in $\Rd.$ Consider the initial
mass distribution $\mu$ as fixed and deterministic and choose the
initial positions $\zni$ in such a way that
\begin{equation*}
    \limn \frac{1}{n}\sum_{i=1}^n\delta_{\zni}=\mu.
\end{equation*}
The empirical measure of the final positions is
\begin{equation*}
   \Nkn=\frac{1}{n} \sum_{i=1}^n \delta_{\Xkni(1)}.
\end{equation*}
It is a random element of $\PR.$ An easy variation of Sanov's
theorem states that for each $k$ the system $\{\Nkn\}_{n\geq 1}$
obeys the LDP in $\PR$ with speed $n$ and the rate function
\begin{equation}\label{eq-11}
   \nu\in\PR\mapsto \inf_{\rho\in\Pi(\mu,\nu)} H(\rho|\pi^k)\in
   [0,\infty].
\end{equation}
Here, $H(\rho|\pi^k)$ is the relative entropy (see (\ref{eq-10})
for its definition) of $\rho$ with respect to $\pi^k$ and
$\pi^k\in\PRR$ is the law of $(Z,Z+Y/\sqrt{k})$ where $Z$ and $Y$
are independent, the law of $Z$ is $\mu$ and $Y$ is a standard
normal vector. On the other hand, $\{Y/\sqrt{k}\}_{k\geq 1}$ obeys
the LDP in $\Rd$ as $k$ tends to infinity with speed $k$ and rate
function $c(u)=|u|^2/2.$
\\
Since
\begin{itemize}
    \item[(i)] the speed of the LDP for $\{Y/\sqrt{k}\}_{k\geq 1}$
    is $k$ and
    \item[(ii)] the rate functions (\ref{eq-11}) and $c(u)=|u|^2/2$ are
reminiscent of $\mathcal{T}_c$ given at (\ref{eq-12}),
\end{itemize}
it wouldn't be surprising that
\begin{itemize}
    \item[(i)] the order of magnitude of $H(\rho|\pi^k)$ is $k$ and
    \item[(ii)] one should mix together two types of LDPs with $n$ and $k$ tending to
    infinity, in order to obtain some LDP with the rate function $\nu\mapsto\mathcal{T}_c(\mu,\nu).$
\end{itemize}
Indeed, denoting for each $\nu\in\PR$ with fixed $\mu,$
\begin{eqnarray*}
  T_k(\nu)&=&\inf_{\rho\in\Pi(\mu,\nu)} H(\rho|\pi^k)/k\quad
\textrm{and} \\
  T(\nu) &=& \mathcal{T}_c(\mu,\nu),
\end{eqnarray*} it will be proved that
the following $\Gamma$-convergence result
\begin{equation}\label{eq-14}
    \Glimk T_k=T
\end{equation}
holds. As a consequence of this convergence result, for each
$\nu\in\PR,$ there exists a sequence $(\nu_k)_{k\geq 1}$ such that
\begin{equation}\label{eq-18}
   \limk \nu_k=\nu\quad \textrm{and} \quad \limk \inf_{\rho\in\Pi(\mu,\nu_k)}
    H(\rho|\pi^k)/k=\mathcal{T}_c(\mu,\nu).
\end{equation}

Theorem \ref{res-09} is the main result of the paper. It states
that $\{\Nkn\}_{k,n\geq 1}$ obeys the doubly indexed LDP as $n$
first tends to infinity, then $k$ tends to infinity with speed
$kn$ and rate function $T,$ see Definition \ref{bb} for the notion
of doubly indexed LDP.

\subsection{An approximation procedure}\label{sec-approx}

The $\Gamma$-limit (\ref{eq-14}) suggests that the sequence of
minimizers $\rho^*_k$ of $H(\rho|\pk)$ subject to the constraint
$\rho\in \Pi(\mu,\nu)$ should converge as $k$ tends to infinity to
some minimizer of $\rho\mapsto\int_{\Rdd}c\,d\rho$ subject to the
same constraint $\rho\in \Pi(\mu,\nu).$ This fails in many
situations. Consider for instance a purely atomic initial measure
$\mu$ and a family of atomic probability measures $\pk.$ Although
$T(\nu)$ may be finite for some diffuse final measure $\nu,$ there
are no $\rho$ in $\Pi(\mu,\nu)$ which are absolutely continuous
with respect to $\pk$ since $\pk_1$ is atomic. Hence,
$T_k(\nu)=+\infty$ for all $k,$ and there are no minimizers
$\rho^*_k$ at all. To take this phenomenon into account, one can
think of the minimization problems
\begin{equation}\label{eq-MKk}
     \textrm{minimize }\quad H(\rho|\pk)/k
     \quad\textrm{subject to}\quad \rho\in \Pi(\mu,\nu_k)\tag{MK$_k$}
\end{equation}
where $(\nu_k)_{k\geq 1}$ satisfies (\ref{eq-18}). I didn't
succeed in proving that $\limk
(\hbox{\ref{eq-MKk}})=(\hbox{\ref{eq-MK}})$ in the sense of
$\Gamma$-convergence.

Alternately, one can relax the constraint $\rho_1=\nu$ by means of
a continuous penalization sequence and consider the three
minimization problems
\begin{align}
    &\textrm{minimize\quad} H(\rho|\pk)/k+\alpha d(\rho_1,\nu) \textrm{\quad subject
    to\quad} \rho_0=\mu\tag{MK$_k^\alpha$}\label{eq-MKka}\\
    &\textrm{minimize\quad} \int_{\Rdd}c\,d\rho+\alpha d(\rho_1,\nu) \textrm{\quad subject
    to\quad} \rho_0=\mu\tag{MK$^\alpha$}\label{eq-MKa}\\
    &\textrm{minimize\quad} \int_{\Rdd}c\,d\rho \textrm{\quad subject
    to\quad} \rho\in\Pi(\mu,\nu)\tag{MK}
\end{align}
where $k,\alpha\geq 1$ are intended to tend to infinity and
$d(\rho_1,\nu)$ is some distance between $\rho_1$ and $\nu$ which
is compatible with the narrow topology of $\PR.$

Note that (\ref{eq-MKka}) is a strictly convex problem while
(\ref{MKa}) and (\ref{MK}) are not. As a consequence
(\ref{eq-MKka}) admits a unique minimizer $\rho^\alpha_k$ while
(\ref{MKa}) and (\ref{MK}) may admit several ones. It will proved
by means of another $\Gamma$-convergence result that
\begin{equation}\label{eq-15}
    \limk (\hbox{\ref{eq-MKka}})=(\hbox{\ref{eq-MKa}}) \quad \textrm{and} \quad \lima (\hbox{\ref{eq-MKa}})=(\hbox{\ref{eq-MK}}).
\end{equation}
These  formulas are to be understood at a formal level. It means
in particular that for each $\alpha,$ $\limk\inf
(\hbox{\ref{eq-MKka}})=\inf (\hbox{\ref{eq-MKa}})$ and all the
limit points of the relatively compact sequence
$(\rho^\alpha_k)_{k\geq 1}$ are minimizers of the limiting
 problem (\ref{eq-MKa}). Similarly, $\lima\inf
(\hbox{\ref{eq-MKa}})=\inf (\hbox{\ref{eq-MK}})$ and denoting
$\rho^\alpha$ a minimizer of (\ref{eq-MKa}),  any limit point of
the relatively compact sequence $(\rho^\alpha)_{\alpha\geq 1}$ is
a minimizer of the limiting  problem (\ref{eq-MK}).

\subsection{Some comment about the results of this paper}
The doubly indexed LDP for  $\{\Nkn\}_{k,n\geq 1},$ the limit
(\ref{eq-14}) and the approximation procedure (\ref{eq-15}) are
new results. Large deviations have only been used as a guideline
to obtain the analytical results (\ref{eq-14}) and (\ref{eq-15}).

In the rest of the paper not only the quadratic cost is considered
but a much wider class of cost functions. In particular, the above
mentioned results hold true for the usual power cost functions
$c(u)=|u|^p$ with $p>0.$ Note that the convexity of $c$ is not
required.

We choosed $\Rd$ as the surrounding space to make the presentation
of the results easier. It is by no way a limitation. Our main
large deviation result (Theorem \ref{res-33}) is stated with
Polish spaces. On the other hand, the proofs of our convergence
results mainly rely on $\Gamma$-convergence. We have done them in
$\Rd,$ but their extension to a Polish space is obvious.

As a by-product of our approach, the Kantorovich duality
(\cite{Vill03}, Theorem 1.3) is recovered, see Theorems
\ref{res-33} and \ref{res-10}. This provides a new proof of it,
although not the shortest one.

\subsection{Literature}
Since Brenier's note \cite{Brenier87} in 1987 which was motivated
by fluid mechanics, optimal transport is a very active area of
applied mathematics. For a comprehensive account on optimal
transport theory, we refer to the monographs of Rachev and
R\"uschendorf \cite{RacRus} and Villani \cite{Vill03}. Villani's
recent Saint-Flour lecture notes \cite{Vill05} are up-to-date and
aimed at a probalistic reader. They introduce newly born
techniques and offer a very long reference list.

Although optimal transport has important consequences in
probability theory (Wasserstein's metrics or transportation
inequalities for instance), it has seldom been studied from a
probabilistic point of view. Let us cite among others the
contributions of Feyel and \"Ust\"unel \cite{FeyUstu04a},
\cite{FeyUstu04b} about the Monge-Kantorovich problem on Wiener
space. Recently, Mikami \cite{Mikami04} has obtained a
probabilistic proof of the existence of a solution to Monge's
problem with a quadratic cost by means of an approximation
procedure by $h$-processes. His approach is based on optimal
control techniques.

Doubly indexed LDs of empirical measures have been studied by
Boucher, Ellis and Turkington in \cite{BET99}. In \cite{Leo05a},
the tight connection between doubly indexed LDs and the
$\Gamma$-convergence of LD rate functions is stressed. This will
be used in the present article.

\subsection{$\Gamma$-convergence}

The $\Gamma$-convergence is a useful tool which is going to be
used repeatedly. We refer to the monograph of G.~Dal~Maso
\cite{DalMaso} for a clear exposition of the subject. Precise
references to the invoked theorems in \cite{DalMaso} will be
written all along the paper.
\\
Recall that if it exists, the $\Gamma$-limit of the sequence
$(f_n)_{n\geq 1}$ of $(-\infty,\infty]$-valued functions on a
topological space $X$ is given for all $x$ in $X$ by
\begin{equation*}
    \Glimn f_n(x)=\sup_{V\in \mathcal{N}(x)}\limn \inf_{y\in V}
    f_n(y)
\end{equation*}
where $\mathcal{N}(x)$ is the set of all neighbourhoods of $x.$
This notion of convergence is well-designed for minimization
problems. Denoting $f=\Glimn f_n$ and taking $(x_n)$  a converging
sequence of minimizers of $(f_n)$ with $\limn x_n=x^*,$ if
$(f_n)_{n\geq 1}$ is equi-coercive we have $\limn \inf f_n=\inf f$
 and  $x^*$ is a minimizer of $f.$

\subsection{Some notations and conventions}\label{sec-conventions}
Let us fix some notations and conventions.
\\
\noindent \textit{Topological conventions}.\
 The space of all
continuous bounded functions on a topological space $\X$ is
denoted by $\CX$ and is equipped with the uniform norm
$\|f\|=\sup_{x\in\X}|f(x)|,$ $f\in\CX.$ Unless specified, its dual
space $\CX'$ is equipped with the $\ast$-weak topology
$\sigma(\CX',\CX).$
\\
Any Polish space $\X$ is equipped with its Borel $\sigma$-field
and  the set $\PX$  of all the probability measures on $\X$  is
equipped with the narrow topology $\sigma(\PX,\CX):$ the relative
topology of $\CX'$ on $\PX.$ While considering random probability
measures, it is necessary to equip $\PX$ with some $\sigma$-field:
we take its Borel $\sigma$-field.

\par\medskip\noindent
\textit{Large deviations.}\ Let $\{V_n\}_{n\geq 1}$ be a sequence
random variables taking their values in some topological space
$\mathcal{V}$ equipped with some $\sigma$-field. One says that
$\{V_n\}_{n\geq 1}$ obeys the Large Deviation Principle (LDP) in
$\mathcal{V}$ with speed $n$ and rate function
$I:\mathcal{V}\to[0,\infty],$ if $I$ is \lsc\ and for all
measurable subset $A$ of $\mathcal{V},$ we have
\begin{eqnarray}\label{eq-ldp}
  -\inf_{v\in \inter A} I(y) &\leq& \lin\nlog (V_n\in A) \\
\nonumber  &\leq& \lsn\nlog (V_n\in A)\leq -\inf_{v\in \cl A}
I(y)\nonumber
\end{eqnarray}
where $\inter A$ and $\cl A$ are the interior and the closure of
$A$ in $\mathcal{V}.$
\\
To emphasize the parameter $n,$ one says that this is a \nldp. If
$\rho_n$ denotes the law of $V_n,$ one also writes that
$\{\rho_n\}_{n\geq 1}$ obeys the \nldp\ in $\mathcal{V}$ with the
rate function $I.$
\\
The rate function $I$ is said to be a good rate function if for
each $a\geq 0,$ the level set $\{I\leq a\}$ is a compact subset of
$\mathcal{V}.$ We shall equivalently write that $I$ is
inf-compact.

\subsection{Organization of the paper}
At Section \ref{sec-main results} the main results are stated
precisely without proof. Their proofs are postponed to Section
\ref{sec-applications}. They rely on preliminary results obtained
at Sections \ref{sec-doubly indexed} and \ref{sec-doubly main}
where general large deviation results are derived for doubly
indexed sequences of random probability measures with our optimal
transport problems in  mind. As a preliminary approach, Section
\ref{sec-simply indexed} is dedicated to easier analogous large
deviation results in terms of simply indexed sequences. Finally,
Section \ref{sec-appendix} is an appendix dedicated to the proof
of a  result about the $\Gamma$-convergence of convex functions
which is used in Section \ref{sec-doubly indexed}. Since we didn't
find this result in the literature, we give its detailed proof.

\tableofcontents

\section{Statement of the results}\label{sec-main results}

The main result of the paper is Theorem  \ref{res-33}, it is
stated in an abstract setting with general Polish spaces. In the
present section, it is restated at Theorem \ref{res-09} without
proof in the particular framework of the optimal transport on
$\Rd.$  All the results of the present section are proved at
Section \ref{sec-applications}, using the results of Sections
\ref{sec-simply indexed}, \ref{sec-doubly indexed} and
\ref{sec-doubly main}.

\subsection{Some transportation cost functions are LD rate functions}

Take a triangular array $(\zni\in\Rd; 1\leq i\leq n, n\geq 1)$ in
$\Rd$ which satisfies
\begin{equation}\label{dd}
    \limn \frac{1}{n}\sum_{i=1}^n\delta_{\zni}=\mu
\end{equation}
for some $\mu\in\PR.$
\\
For each $z\in\Rd,$ let $\{U^k_z\}_{k\geq 1}$ be a sequence of
$\Rd$-valued random variables. For each $k\geq 1$ and $n\geq 1,$
take $n$ \emph{independent} random variables $(\Ykni)_{1\leq i\leq
n}$ where
\begin{equation}\label{eq-13a}
    \Ykni\eqlaw U^k_{\zni}.
\end{equation}
 For each $k,$ $(\Ykni; 1\leq i\leq n)_{n\geq
1}$ is a triangular array of independent particles which, in the
general case, are not identically distributed because of the
contribution of the deterministic $\zni$'s. An important example
is given by $U^k_z=z+U^k$ with $\{U^k\}_{k\geq 1}$  a sequence of
$\Rd$-valued random variables. This gives for each $k,n\geq 1$
\begin{equation}\label{eq-13b}
    \Ykni=\zni + U^k_i
\end{equation}
where $(U^k_i)_{1\leq i\leq n}$ are independent copies of $U^k.$

We are interested in the large deviations of the empirical
measures on $\Rd$
\begin{equation}\label{eq-04}
   \Nkn=\frac{1}{n} \sum_{i=1}^n \delta_{\Ykni}
\end{equation}
as $n$ first tends to infinity, then $k$ tends to infinity. More
precisely, doubly indexed LDPs in the sense of the following
definition will be proved.

\begin{definition}[Doubly indexed LDP]\label{bb}
Let $\PX$ be the set of all probability measures built on the
Borel $\sigma$-field of a Polish space $\X.$ The set $\PX$ is
equipped with the topology of narrow convergence and with the
corresponding Borel $\sigma$-filed.
\\
  One says that a doubly indexed $\PX$-valued sequence $\{\Lkn\}_{k,n\geq 1}$
  obeys the \knldp\ in $\PX$ with the rate function
  $I:\PX\to[0,\infty],$ if for all measurable subset $B$ of $\PX,$ we have
\begin{eqnarray}\label{ee}
  -\inf_{Q\in \inter B} I(Q) &\leq& \likn\knlog (\Lkn\in B)\nonumber \\
  &\leq& \lskn\knlog (\Lkn\in B)\leq -\inf_{Q\in \cl B} I(Q)
\end{eqnarray}
where $\inter B$ and $\cl B$ are the interior and closure of $B$
in $\PX.$
\end{definition}

\begin{assumptions}\label{assumptions}
This set of assumptions holds for the present
section and Section \ref{sec-applications}.
\begin{itemize}
    \item (\ref{dd}) holds for some $\mu$ in $\PR,$
    \item for each $k\geq 1,$ $(Law(U^k_z);z\in\Rd)$ is a Feller
    system in the sense of Definition \ref{def-Feller} below and
    \item for each $z\in\Rd,$ $\{U^k_z\}_{k\geq 1}$ obeys the \kldp\ in $\Rd$ with the
    \emph{good} rate function $c_z(u)\in [0,\infty],$ $u\in\Rd.$
\end{itemize}
\end{assumptions}

\begin{definition}\label{def-Feller}
Let $\Z$ and $\X$ be two topological spaces. The system of Borel
probability measures $(P_z; z\in\Z)$ on $\X$ is a \emph{Feller
system} if for all $f$ in $\CX,$ $z\in\Z\mapsto \IX
f(x)\,P_z(dx)\in \mathbb{R}$ is a continuous function on $\Z.$
\end{definition}

The next theorem shows that the large deviations of $\{\Nkn\}$ are
closely related to optimal transport.

\begin{thm}\label{res-09}
  The doubly indexed system $\{\Nkn\}_{k,n\geq 1}$
obeys the \knldp\ in $\PR$ with the rate function
\begin{equation*}
    T(\nu)=\mathcal{T}_c(\mu,\nu)
\end{equation*}
for all $\nu\in\PR,$ where the cost function is given by
\begin{equation}\label{eq-c}
    c(x_0,x_1)=c_{x_0}(x_1),\quad x_0,x_1\in\Rd.
\end{equation}
In the special case where (\ref{eq-13b}) holds and $\{U^k\}_{k\geq
1}$ obeys the \kldp\ in $\Rd$ with the good rate function
$c:\Rd\to[0,\infty],$ we have $c(x_0,x_1)=c(x_1-x_0),$
$x_0,x_1\in\Rd.$
\end{thm}

\begin{examples}\label{ex-02}
In the special case where (\ref{eq-13b}) holds, we give some
examples of $\{U^k\}$ and the corresponding cost function $c.$
\begin{enumerate}
    \item With $U^k=Y/\sqrt k$ where $Y$ is a standard normal
    random vector on $\Rd,$ we get
$$
c(u)=|u|^2/2,\quad u\in\Rd.
$$
    This
    is the usual quadratic cost function.
    \item Let $(Y_m)_{m\geq 1}$ be a sequence of independent
    copies of a $\Rd$-valued random vector $Y$ which satisfies
    $\E e^{a|Y|}<\infty$ for some $a>0.$
    With $U^k=\frac{1}{k}\sum_{1\leq m\leq k}Y_m,$
    Cram\'er's theorem (\cite{DZ}, Corollary 6.1.6) states that $\{U^k\}$
    obeys the \kldp\ in $\Rd$ with the rate function $c=\cY:$
\begin{equation}\label{cc}
    \cY(u)=\sup_{\zeta\in\Rd}\{\langle \zeta,u\rangle-\log\mathbb{E}e^{\langle\zeta,Y
    \rangle}\},\quad u\in\Rd.
\end{equation}
    Observe that (1) is a specific instance of (2).
    \item Let $(Y_m)_{m\geq 1}$ be as above and let $\alpha$ be
    any continuous mapping on $\Rd.$ With $U^k=\alpha(\frac{1}{k}\sum_{1\leq m\leq
    k}Y_m)$ we obtain $c(u)=\inf\{\cY(v);v\in \Rd, \alpha(v)=u\}, u\in\Rd$
    as a consequence of the contraction principle.
    In particular if $\alpha$ is a continuous injective mapping,
    then
\begin{equation*}
    c=\cY\circ \alpha^{-1}.
\end{equation*}
    \item For instance, mixing (1) and (3) with $\alpha=\alpha_p$
    given for each $p>0$ and $v\in\Rd$ by
    $\alpha_p(v)=2^{-1/p}|v|^{2/p-1}v,$ taking
    $U^k=(2k)^{-1/p}|Y|^{2/p-1}Y$ where $Y$ is a standard normal
    random vector on $\Rd,$ we get
$$
c(u)=|u|^p,\quad u\in\Rd.
$$
    Note that $U^k\eqlaw k^{-1/p}Y_p$ where the density of the law
    of $Y_p$ is $C|z|^{p/2-1}e^{-|z|^p}.$
\end{enumerate}
\end{examples}

\begin{examples}\label{ex-03}
We recall some well-known examples of Cram\'er transform $\cY.$
\begin{enumerate}
    \item To obtain the quadratic cost function $\cY(u)=|u|^2/2,$ choose $Y$
    as a standard normal random vector in $\Rd.$
    \item Taking $Y$ such that
    $\mathbb{P}(Y=+1)=\mathbb{P}(Y=-1)=1/2,$ leads to\\
    $\cY(u)=\left\{%
\begin{array}{ll}
    [(1+u)\log(1+u)+(1-u)\log(1-u)]/2, & \hbox{if }-1<u<+1 \\
    \log 2, & \hbox{if } u\in\{-1,+1\}\\
    +\infty, & \hbox{if } u\not\in [-1,+1]. \\
\end{array}%
\right.    $
    \item If $Y$ has an exponential law with expestation 1, $\cY(u)=u-1-\log
    u$ if $u>0$ and $\cY(u)=+\infty$ if $u\leq 0.$
    \item If $Y$ has a Poisson law with expectation 1, $\cY(u)=u\log
    u-u+1$ if $u>0,$ $\cY(0)=1$ and $\cY(u)=+\infty$ if $u< 0.$
     \item We have $\cY(0)=0$ if and only if $\mathbb{E}Y=0.$
     \item More generally, $\cY(u)\in [0,+\infty]$ and $\cY(u)=0$ if
     and only if $u=\mathbb{E}Y.$
     \item We have
     $c^{aY+b}(u)=c^Y((u-b)/a)$ for all real $a\not=0$ and $b\in\Rd.$
\end{enumerate}
\end{examples}

\begin{examples}\label{ex-01}
 If $\mathbb{E}Y=0,$ $\cY$ is quadratic at the origin since $\cY(u)=\langle
 u,\Gamma_Y^{-1}
    u\rangle/2+o(|u|^2)$ where $\Gamma_Y$ is the covariance of
    $Y.$ This rules out the usual costs $c(u)=|u|^p$ with $p\not=2.$
    \\
    Nevertheless, taking $Y$ a real valued variable with  density
    $C\exp(-|z|^p/p)$ with $p\geq 1$ leads to
    $\cY(u)=|u|^p/p(1+o_{|u|\rightarrow\infty}(1)).$
    The case $p=1$ follows from Example \ref{ex-03}-(3) above. To see that the result still holds with $p>1,$
    compute by means of the Laplace method the principal part as $\zeta$ tends to infinity of $\int_0^\infty
    e^{-z^p/p}e^{\zeta z}\,dz=\sqrt{2\pi(q-1)}\zeta^{1-q/2}e^{\zeta^q/q}(1+o_{\zeta\rightarrow
    +\infty}(1))$ where $1/p+1/q=1.$
    \\
    Of course, we deduce a related $d$-dimensional result considering
    $Y$ with the density
    $C\exp(-|z|_p^p/p)$ where $|z|_p^p=\sum_{i\leq d}|z_i|^p.$ This gives
    $\cY(u)=|u|_p^p/p(1+o_{|u|\rightarrow\infty}(1)).$
\end{examples}

The drawback of the specific shape of any Cram\'er's transform
$\cY$ (see Examples \ref{ex-01}) is overcome by means of a
continuous transformation as in Examples \ref{ex-02}-(3 \& 4).

\subsection{Convergence results}
The structure of (\ref{ee}) suggests that a \knldp\ may be seen as
the limit as $k$ tends to infinity of \nldp s indexed by $k.$ This
is true and made precise at Proposition \ref{res-03a} and Theorem
\ref{res-04a} below.
\\
Let us have a look at the \nldp\ satisfied by $\{\Nkn\}_{n\geq 1}$
with $k$ fixed. It is very similar to the \nldp\  of Sanov's
theorem, see Proposition \ref{res-03a} below. The only difference
comes from the contribution of the initial positions $\zni$ which
make $(\Ykni)$ a triangular array of \emph{non}-identically
independent variables. Recall that Sanov's theorem (\cite{DZ},
Theorem 6.2.10) states that the empirical measures
$\{\frac{1}{n}\sum_{i=1}^n \delta_{X_i}\}_{n\geq 1}$ of a sequence
of independent $P$-distributed random variables taking their
values in a Polish space $\X$ obey the \nldp\ in $\PX$ with the
rate function
\begin{equation}\label{eq-10}
  Q\in\PX\mapsto  H(Q|P)=\left\{%
\begin{array}{ll}
    \int_\X \log\left(\frac{dQ}{dP}\right)\,dQ & \hbox{if $Q\prec P$} \\
    +\infty & \hbox{otherwise.} \\
\end{array}%
\right.
\end{equation}
$H(Q|P)$ is called the relative entropy of $Q$ with respect to
$P.$

Consider now the random empirical measures on
$\Rdd=\mathbb{R}^d\!\times\!\mathbb{R}^d$ which are defined by
\begin{equation}\label{eq-03}
    \Mkn=\frac{1}{n}\sum_{i=1}^n \delta_{(\zni,\Ykni)}
\end{equation}
for all $k,n\geq 1.$ Clearly, $\Nkn$ is the second marginal of
$\Mkn.$  Denote for each $k\geq 1$
\begin{equation}\label{eq-02}
    \pk(dx_0dx_1)=\mu(dx_0)Law(U^k_{x_0})(dx_1)\in\PRR
\end{equation}
This means that $\pk=Law(X(0),X^k(1))$ where $X(0)$ is a
$\Rd$-valued random variable which is $\mu$-distributed and
$\P(X^k(1)\in dx_1\mid X(0)=x_0)=Law(U^k_{x_0})(dx_1)$. Define
\begin{equation*}
    S_k(\rho)=\left\{%
\begin{array}{ll}
   \frac{1}{k} H(\rho|\pk) & \hbox{if }\rho_0=\mu \\
    +\infty & \hbox{otherwise} \\
\end{array}%
\right.,\quad \rho\in\PRR
\end{equation*}
and
\begin{equation}\label{a1a}
   T_k(\nu)= \inf_{\rho\in \Pi(\mu,\nu)}\frac{1}{k}H(\rho|\pk), \quad
   \nu\in\PR.
\end{equation}

\begin{prop}\label{res-03a}
For each fixed $k\geq 1,$
\begin{enumerate}
    \item[(a)]
    $\{\Mkn\}_{n\geq 1}$ obeys the \nldp\ in $\PRR$ with the good rate
    function $kS_k$ and
    \item[(b)]
    $\{N^k_n\}_{n\geq 1}$ obeys the \nldp\ in $\PR$ with the good rate
    function $kT_k.$
\end{enumerate}
\end{prop}

The order of magnitude of $H(\rho|\pk)$ is $k,$ since $\{U^k\}$
obeys a \kldp. The rescaled entropy $S_k$ is of order 1. If it
exists, $\limk S_k$ may be interpreted as a specific entropy (see
\cite{Rue78}). It happens that $S_k$ and $T_k$ $\Gamma$-converge.
The limit of $S_k$ is
\begin{equation*}
    S(\rho)=\left\{%
\begin{array}{ll}
   \int_{\Rdd}c\,d\rho & \hbox{if }\rho_0=\mu \\
    +\infty & \hbox{otherwise} \\
\end{array}%
\right.,\quad \rho\in\PRR
\end{equation*}
where $\int_{\Rdd}c\,d\rho=\int_{\Rdd} c(x_0,x_1)\,\rho(dx_0dx_1)$
and $c$ is given at (\ref{eq-c}).

\begin{thm}\label{res-04a}
We have
\begin{enumerate}
    \item[(a)] $\Glimk S_k=S$ in $\PRR$ and
    \item[(b)] $\Glimk T_k=T$ in $\PR.$
\end{enumerate}
\end{thm}
These limits will allow us to deduce the following approximation
results. Recall that the minimization problems (\ref{MKka}),
(\ref{MKa}) and (\ref{MK}) are defined at Section
\ref{sec-approx}.

\begin{thm}\label{res-14a}
Assume that $\mathcal{T}_c(\mu,\nu)<\infty.$
\begin{itemize}
    \item[(a)] We have:\quad
$\lima\limk
\inf_{\rho\in\Pi_0(\mu)}\left\{\frac{1}{k}H(\rho|\pk)+\alpha
d(\rho_1,\nu)\right\}=\mathcal{T}_c(\mu,\nu).$
    \item[(b)] For each $k$ and $\alpha,$ (\ref{MKka}) admits a unique solution
    $\rho^\alpha_k$ in $\PRR.$ For each $\alpha,$ $(\rho_k^\alpha)_{k\geq
    1}$ is a relatively compact sequence in $\PRR$ and  any limit point of $(\rho_k^\alpha)_{k\geq
    1}$ is a solution of (\ref{MKa}).
    \item[(c)] For each $\alpha,$ (\ref{MKa}) admits at least a (possibly not unique) solution $\rho^\alpha.$
    The sequence $(\rho^\alpha)_{\alpha\geq 1}$ is relatively compact in
    $\PRR$ and  any limit point of $(\rho^\alpha)_{\alpha\geq
    1}$ is a solution of (\ref{MK}).
\end{itemize}
\end{thm}

\subsection{The proofs} The proofs of these announced results are
done at Section \ref{sec-applications}. Theorem \ref{res-09} is
part of Theorem \ref{res-10}, Proposition \ref{res-03a} is Lemma
\ref{res-12}, Theorem \ref{res-04a} is Theorem \ref{res-13} and
Theorem \ref{res-14a} is Theorem \ref{res-14}.

\section{Large deviations of a simply indexed sequence of random measures}
\label{sec-simply indexed}

As a warming-up exercice, let us first consider a usual sequence
of random measures.

We present an abstract setting instead of the situation described
at Section \ref{sec-main results}. Let $\X$ and $\Z$ be two Polish
spaces which play respectively the part of the space of "paths"
$\Rdd$ and the space of initial conditions $\Rd.$  The cost of
this extension is quite low: the main property of Polish spaces to
be used later is that any Borel probability measure is tight.

Take a triangular array
 $(\zni\in\Z; 1\leq i\leq n, n\geq 1)$
on $\Z$ such that the sequence of empirical measures
$\mu_n=\frac{1}{n}\sum_{i=1}^n \delta_\zni \in\PZ$ satisfies
\begin{equation}\label{eq-mun}
\limn \mu_n=\mu
\end{equation}
for some probability measure $\mu\in\PZ.$ Let $(P_z\in\PX;
z\in\Z)$ be  a collection of probability laws on $\X$ which is
assumed to be a Feller system in the sense of Definition
\ref{def-Feller}.
\\
We work  with a triangular array of \emph{independent} $\X$-valued
random variables $(\Xni; 1\leq i\leq n, n\geq 1)$ where for each
index $(n,i)$ the law of $\Xni$ is $P_\zni.$ This means that for
all $n\geq 1,$
\begin{equation*}
    \mathcal{L}aw(\Xni; 1\leq i\leq n)=\otimes_{1\leq i\leq n}
    P_\zni.
\end{equation*}
Proposition \ref{res-22} below states a LDP in $\PX$ for the
empirical measures
\begin{equation*}
    L_n=\frac{1}{n} \sum_{i=1}^n \delta_\Xni
\end{equation*}
as $n$ tends to infinity. It is a variant of Sanov's theorem which
has already been studied by Dawson and G\"artner in \cite{DG87}
and revisited by Cattiaux and L\'eonard in \cite{CL95}.
Nevertheless, the expression (\ref{eq-21}) of the rate function
doesn't appear in these cited papers.  The proof of Proposition
\ref{res-22} will be done as a first step for the proof of the LDP
of a doubly indexed sequence: most of its ingredients will be
recycled at Section \ref{sec-doubly indexed}.

 \par\medskip\noindent \textbf{Notations.}\
 We write shortly $\PZX$ and
$\CZX$ for $\mathcal{P}_{\ZX}$ and $C_{\ZX}.$ The dual space
$\CZX'$ of $(\CZX,\|\cdot\|)$ is equipped with the $\ast$-weak
topology $\sigma(\CZX',\CZX),$ see Section \ref{sec-conventions}.
\\
Let $(Z,X)$ be the canonical projections: $Z(z,x)=z, X(z,x)=x,$
$(z,x)\in\ZX.$ For any  $q\in\PZX,$ we write the desintegration
$$
q(dzdx)=q_\Z(dz)q^z(dx)
$$
where $q_\Z(dz)=q(Z\in dz)$ is the (marginal) law of $Z$ under $q$
and $q^z(dx)=q(X\in dx\mid Z=z), z\in\Z,$ is a regular conditional
version of the law of $X$ knowing that $Z=z.$ We also define
$p\in\PZX$ by
$$
p(dzdx) = \mu(dz)P_z(dx).
$$

The LDP for $\{L_n\}_{n\geq 1}$ will be obtained as a direct
consequence of the contraction principle applied to some LDP for
the sequence of $\PZX$-valued random variables
\begin{equation*}
    K_n=\frac{1}{n} \sum_{i=1}^n \delta_{(\zni,\Xni)}, \quad n\geq
    1.
\end{equation*}

\begin{prop}\label{res-21}
Suppose that (\ref{eq-mun}) holds for some $\mu$ in $\PZ$ and that
$(P_z; z\in\Z)$ is a Feller system. Then $\{K_n\}_{n\geq 1}$ obeys
the LDP in $\PZX$ with the good rate function
\begin{equation}\label{eq-22}
h(q):=\left\{%
\begin{array}{ll}
    H(q|p)=\IZ H(q^z|P_z)\,\mu(dz) & \hbox{if }q_\Z=\mu \\
    +\infty & \hbox{otherwise} \\
\end{array}%
\right.,\quad q\in\PZX
\end{equation}
\end{prop}

 \proof
For all $n$ and all $F\in \CZX,$ the normalized log-Laplace
transform of $K_n$ is
\begin{eqnarray*}
    \psi_n(F)
  &:=&\un \log \E \exp(n\langle F,K_n\rangle)\\
  &=& \un\sum_{i=1}^n \log \E e^{F(\zni,\Xni)} \\
  &=& \IZ \log\langle e^{F_z},P_z\rangle\, \mu_n(dz).
\end{eqnarray*}
As $(\mu_n)_{n\geq 1}$ converges to $\mu$ and $(P_z;z\in\Z)$ is a
Feller system, for all $F\in \CZX$ we have the limit:
\begin{eqnarray}\label{eq-psi}
  \psi(F) &:=& \limn \psi_n(F)\nonumber \\
   &=& \IZ \log\langle e^{F_z},P_z\rangle\, \mu(dz).
\end{eqnarray}

Following the proof of Sanov's theorem (see \cite{DZ}, Section
6.4) based on Dawson-G\"artner's theorem on the projective limit
of LD systems (see \cite{DG87}, Section 3), one obtains that
$\{K_n\}$ obeys the LDP in $\CZX'$ with the rate function

\begin{equation}\label{eq-203}
    \psi^*(q)=\sup_{F\in \CZX}\left\{\langle F,q\rangle -\IZ
\log\langle e^{F_z},P_z\rangle\, \mu(dz)\right\}, \quad q\in
\CZX'.
\end{equation}

It is proved at Lemma \ref{res-202} below, that for all $q\in
\CZX',$
\begin{equation*}
  \psi^*(q):=\left\{%
\begin{array}{ll}
    h(q) & \hbox{if }q\in\PZX \\
    +\infty & \hbox{otherwise} \\
\end{array}%
\right.
\end{equation*}

It follows that $\{K_n\}_{n\geq 1}$ obeys the LDP in $\PZX$ with
the rate function $h.$

It remains to note that as the relative entropy is inf-compact and
$\{q\in\PZX; q_\Z=\mu\}$ is closed, $h$ is also inf-compact: it is
a good rate function.
\endproof

As a by-product of this proof, we have the following corollary
which is mentioned for future use.

\begin{cor}\label{res-201} [Hypotheses of Proposition
\ref{res-21}].  The random system $\{K_n\}_{n\geq 1}$ obeys the
LDP in $\CZX'$ with the rate function $\psi^*$ given at
(\ref{eq-203}).
\end{cor}

During the proof of Proposition \ref{res-21} we have used the
following lemma.

\begin{lemma}\label{res-202}
With $\psi^*(q)$ defined by formula (\ref{eq-203})  we have
\begin{equation*}
\psi^*(q)=\left\{%
\begin{array}{ll}
    H(q|p)=\IZ H(q^z|P_z)\,\mu(dz) & \hbox{if }q\in\PZX \hbox{  and  } q_\Z=\mu \\
    +\infty & \hbox{otherwise.} \\
\end{array}%
\right.
\end{equation*}
\end{lemma}

\proof  The proof is twofold. We show that
\begin{itemize}
    \item[(i)] for all $q\in \CZX',$ $\psi^*(q)<+\infty$
    implies that $q$ belongs to $\PZX$ and its $z$-marginal is
    $q_\Z=\mu.$
    \item[(ii)] Then, we show that for all $q\in\PZX$ such that
    $q_\Z=\mu,$ we have $\psi^*(q)=H(q|p).$
\end{itemize}

Let $q\in \CZX'$ be such that
\begin{equation*}
    \sup_{F\in \CZX}\{\langle F,q\rangle-\psi(F)\}=\psi^*(q)<\infty.
\end{equation*}

Let us begin with the proof of (i).
\\
 $\bullet$ Let us show that
$q\geq 0.$ Let $F_o\in \CZX$ be such that $F_o\geq 0.$ As $\psi(
aF_o)\leq 0$ for all $a\leq 0,$
\begin{eqnarray*}
  \psi^*(q)
  &\geq& \sup_{a\leq 0}\{a\langle F_o,q\rangle-\psi( aF_o)\} \\
  &\geq& \sup_{a\leq 0}\{a\langle F_o,q\rangle\}\\
  &=& \left\{%
\begin{array}{ll}
    0, & \hbox{if } \langle F_o,q\rangle\geq 0\\
    +\infty, & \hbox{otherwise} \\
\end{array}%
\right.
\end{eqnarray*}
Therefore, as $\psi^*(q)<\infty,$ $\langle F_o,q\rangle\geq 0$ for
all $F_o\geq 0,$ which is the desired result.
\\
$\bullet$  Let us show that $\langle \1,q\rangle=1.$ For any
constant function $F\equiv c\in \mathbb{R},$ we have
$\psi(c\1)=c.$ It follows that
\begin{eqnarray*}
  \psi^*(q)
  &\geq& \sup_{c\in\mathbb{R}}\{c\langle \1,q\rangle-\psi(c\1)\} \\
  &\geq& \sup_{c\in\mathbb{R}}\{c(\langle \1,q\rangle-1)\}\\
  &=& \left\{%
\begin{array}{ll}
    0, & \hbox{if } \langle \1,q\rangle=1\\
    +\infty, & \hbox{otherwise}\\
\end{array}%
\right.
\end{eqnarray*}
from which the result follows.
\\
$\bullet$ In order to prove that $q$ is $\sigma$-additive, we
 have to prove that for any sequence $(F_n)_{n\geq 1}$ in
$\CZX$ such that $F_n\geq 0$ for all $n$ and $(F_n(z,x))_{n\geq
1}$ decreases to zero for each $(z,x)\in\ZX,$ we have
\begin{equation}\label{eq-204}
    \lim_{n\rightarrow\infty} \langle F_n,q\rangle=0.
\end{equation}
For such a sequence, one can apply the dominated convergence
theorem to obtain that
\[
\lim_{n\rightarrow\infty}\psi(aF_n)=0,
\]
for all $a\geq 0.$ It follows that for all $q\in \CZX',$
\begin{eqnarray*}
  \psi^*(q)
  &\geq& \sup_{a\geq 0}\limsup_{n\rightarrow\infty}\{a\langle F_n,q\rangle-\psi(aF_n)\} \\
   &\geq& \sup_{a\geq 0}\left(\limsup_{n\rightarrow\infty} a\langle F_n,q\rangle-\lim_{n\rightarrow\infty}\psi(aF_n)\right) \\
   &=& \sup_{a\geq 0}a\limsup_{n\rightarrow\infty} \langle F_n,q\rangle \\
  &=& \left\{%
\begin{array}{ll}
    0 & \hbox{if }\limsup_{n\rightarrow\infty} \langle F_n,q\rangle\leq 0 \\
    +\infty & \hbox{otherwise.} \\
\end{array}%
\right.
\end{eqnarray*}
Therefore, as $\psi^*(q)<\infty,$ we have
$\limsup_{n\rightarrow\infty} \langle F_n,q\rangle\leq 0.$ Since
we have just seen that $q\geq 0,$ we have obtained (\ref{eq-204}).

This completes the proof of $q\in\PZX$ since we have proved that
any $q\in \CZX'$ such that $\psi^*(q)<\infty$ is nonnegative, has
a unit mass and satisfies (\ref{eq-204}). Therefore, $q$ is
uniquely identified with a probability measure on the Polish space
$\ZX$ (see \cite{Neveu}, Proposition II-7-2).

\par\medskip\noindent
$\bullet$ To complete the proof of (i), it remains to show that
for any $q\in\PZX,$ $\psi^*(q)<\infty$ implies that $q_\Z=\mu.$
Indeed, choosing $F(z,x)=g(z)$ not depending on $x$ with
$g\in\CZ,$ one sees that
\begin{eqnarray*}
  \psi^*(q) &\geq& \sup_{g\in\CZ}\{\langle g, q_\Z\rangle-\langle g,\mu \rangle\} \\
   &=& \left\{%
\begin{array}{ll}
    0, & \hbox{if }q_\Z=\mu \\
    +\infty, & \hbox{otherwise} \\
\end{array}%
\right.
\end{eqnarray*}
which gives the announced result.

\par\medskip\noindent
Now, let us show (ii). For all $q\in\PZX$ such that $q_\Z=\mu$ or
equivalently  such that $q(dzdx)=\mu(dz)q^z(dx),$ we have
\begin{eqnarray}\label{eq-200}
  \psi^*(q)
  &=& \sup_{F\in \CZX}\IZ (\langle F_z,q^z\rangle-\log\langle e^{F_z},P_z\rangle)\,\mu(dz)\nonumber \\
  &\leq& \IZ \sup_{f\in\CX}\{\langle f,q^z\rangle-\log\langle e^{f},P_z\rangle\}\,\mu(dz)\nonumber \\
  &\stackrel{(a)}{=}& \IZ H(q^z|P_z)\,\mu(dz) \nonumber\\
  &\stackrel{(b)}{=}& H(q|p)
\end{eqnarray}
where equality (a) follows from the well-known variational
representation of the relative entropy in a Polish space $\X$
\begin{equation}\label{eq-var rep of H}
    H(Q|P)=\sup_{f\in\CX}\left\{\IX f\, dQ-\log\IX e^f\,
    dP\right\},\quad P,Q\in\PX
\end{equation}
and equality (b) follows from the tensorization property
\begin{equation}\label{eq-tenso}
    H(q|p)=H(q_\Z|p_\Z) + \IZ H(q^z|p^z)\,q_\Z(dz)
\end{equation}
since $p_\Z=q_\Z=\mu$ and $p^z=p(\cdot\mid Z=z)=P_z.$ Note that
$z\mapsto H(q^z|P_z)$ is measurable. Indeed, $(Q,P)\mapsto H(Q|P)$
is measurable as a \lsc\ function and $z\mapsto (q^z,P_z)$ is
measurable since its coordinates are measurable: $z\mapsto q^z$ is
measurable as a regular conditional version in a Polish space and
$z\mapsto P_z$ is assumed to be continuous. We have just proved
that $\psi^*\leq h.$

The converse inequality follows from Jensen's inequality:
$\psi(F)\leq \log \IZ \langle e^{F_z},P_z\rangle\, \mu(dz)=\log
\langle e^F,p\rangle$ for all $F\in \CZX.$ Indeed, taking the
convex conjugates leads us for all $q\in\PZX$ to
\begin{eqnarray}\label{eq-201}
  \psi^*(q)
  &\geq& \sup_{F\in \CZX}\{\langle F,q\rangle-\log \langle e^F,p\rangle\}\nonumber \\
  &=& H(q|p)
\end{eqnarray}
This equality is (\ref{eq-var rep of H}). This completes the proof
of the lemma.
\endproof

\begin{remark}\label{rem-01}
The $\|\cdot\|$-continuity of $q$ in $\CZX'$ didn't play any role
in the proof. Only its linearity has been used.
\end{remark}

Now, we investigate the large deviations of
\begin{equation*}
    L_n=\frac{1}{n} \sum_{i=1}^n \delta_\Xni.
\end{equation*}
Let us denote the $\X$-marginal of $p$ by
\[
P(dx)= \IZ P_z(dx)\,\mu(dz)\in\PX.
\]
\begin{prop}\label{res-22}
Suppose that (\ref{eq-mun}) holds for some $\mu$ in $\PZ$ and that
$(P_z; z\in\Z)$ is a Feller system.
\begin{enumerate}
 \item[(a)]
    Then, $\{L_n\}_{n\geq 1}$ obeys the LDP in $\PX$ with the good rate
    function $H$ which is defined for all $Q\in\PX$ by
    \begin{equation}\label{eq-21}
    H(Q)= \inf\left\{\IZ H(\Pi_z|P_z)\,\mu(dz); (\Pi_z)_{z\in\Z} : \IZ \Pi_z\,\mu(dz)=Q\right\}
    \end{equation}
    where the transition kernels $z\in\Z\mapsto \Pi_z\in\PX$ are measurable.
 \item[(b)]
    If $H(Q)<+\infty,$ there exists a unique (up to $\mu$-a.e. equality) kernel $(\Pi^*_z)_{z\in\Z}$
    which realizes the infimum in (\ref{eq-21}): $H(Q)=\IZ H(\Pi_z^*|P_z)\,\mu(dz).$
 \item[(c)]
    If in addition the Feller system $(P_z)_{z\in\Z}$ satisfies
    \[
    P_z=P (\cdot\mid \beta(X)=z)
    \]
    for $\mu$-almost every $z\in\Z$ and some continuous function
    $\beta: \X\to\Z,$ we have for all $Q\in\PX$
    \begin{equation*}
    H(Q)=\left\{%
    \begin{array}{ll}
    H(Q|P) & \hbox{if $Q$ satifies }\beta\push Q=\mu \\
    +\infty & \hbox{otherwise.} \\
    \end{array}%
    \right.
    \end{equation*}
    and the minimizing kernel $(\Pi^*_z)$ of (\ref{eq-21}) is $\Pi^*_z=Q(\cdot\mid
    \beta(X)=z),$ for $\mu$-almost every $z.$
\end{enumerate}
\end{prop}

\proof

Let us prove (a). As $L_n$ is the $\X$-marginal of $K_n$ and
$\{K_n\}$ obeys the LDP with a good rate function, the statement
(a) follows from the contraction principle (see \cite{DZ}, Theorem
4.2.1): $\{L_n\}$ obeys the LDP in $\PX$ with the good rate
function $H(Q)=\inf\{h(q); q\in\PZX: q_\X=Q\}$ which is
(\ref{eq-21}).

\par\medskip\noindent
 The statement (b) immediately follows from the strict convexity and the
inf-compactness of $q\mapsto H(q |p)$ which is restricted to the
closed convex set $\{q\in\PZX: q_\Z=\mu, q_\X=Q\}.$

\par\medskip\noindent
Let us prove (c). To do this, we rewrite the proof of Proposition
\ref{res-21} with $L_n$ instead of $K_n.$ We obtain that $\{L_n\}$
obeys the LDP in $\PX$ with the rate function
\begin{equation}\label{eq-202}
  \Psi^*(Q)=\sup_{f\in \CX}\left\{\langle f,Q\rangle -\IZ
\log\langle e^{f},P_z\rangle\, \mu(dz)\right\}
\end{equation}
This equality is (\ref{eq-203}) where we replace $q$ by $Q$ and
$F(z,x)$ by $f(x).$ Choosing $f\in \CX$ of the form $f=g\circ
\beta$ with $g$ in $\CZ$ gives us
 for all $Q\in\PX$
\begin{eqnarray*}
  \Psi^*(Q) &\geq& \sup_{g\in \CZ}\{\langle g, \beta\push Q\rangle-\langle g,\mu \rangle\} \\
   &=& \left\{%
\begin{array}{ll}
    0 & \hbox{if }\beta\push Q=\mu \\
    +\infty & \hbox{otherwise} \\
\end{array}%
\right.
\end{eqnarray*}
It follows that $\Psi^*(Q)<\infty$ implies that $\beta\push
Q=\mu.$ For such a $Q,$ as in the proof of inequality
(\ref{eq-200}), we obtain the inequality in  $\Psi^*(Q)\leq \IZ
H\big(Q(\cdot| Z=z)\mid P_z\big)\,\mu(dz)=H(Q|P).$  This last
equality follows from the tensorization property of the relative
entropy, see (\ref{eq-tenso}). This proves that $\Psi^*\leq H.$
The converse inequality follows from Jensen's inequality exactly
as in the proof of inequality (\ref{eq-201}). We have shown that
\begin{equation}\label{eq-210}
    \Psi^*=H.
\end{equation}

The last statement about the minimizing kernel is a direct
consequence of the tensorization formula (\ref{eq-tenso}):
\begin{eqnarray*}
 && \inf\left\{\IZ H(\Pi_z|P_z)\,\mu(dz); (\Pi_z)_{z\in\Z} : \IZ
 \Pi_z\,\mu(dz)=Q\right\}\\
  &=& H(Q| P)\\
  &=& H(Q_\Z | P_\Z)+\IZ H\big(Q(\cdot|Z=z)\mid P(\cdot|Z=z)\big)\,Q_\Z(dz) \\
  &=& \IZ H\big(Q(\cdot|Z=z)\mid P_z\big)\,\mu(dz)
\end{eqnarray*}
where the first equality follows from (a) and the first part of
this statement, and the last equality follows from $H(Q_\Z |
P_\Z)=H(\mu |\mu)=0.$
\endproof

\begin{remark}
The identity (\ref{eq-21}) is a formal inf-convolution formula and
(\ref{eq-202}) is its dual formulation: ``the convex conjugate of
an inf-convolution is the sum of the convex conjugates".
\end{remark}

\begin{remark}\label{rem-21}
Statement (c) holds true also when $\beta$ is only assumed to be
measurable. Indeed, (\ref{eq-202}) can be strengthened by
 \begin{equation*}
  \Psi^*(Q)=\sup_{f\in \CX}\left\{\langle f,Q\rangle -\IZ
\log\langle e^{f},P_z\rangle\, \mu(dz)\right\} =\sup_{f\in
B(\X)}\left\{\langle f,Q\rangle -\IZ \log\langle
e^{f},P_z\rangle\, \mu(dz)\right\},
\end{equation*}
for all $ Q\in\PX,$ where $B(\X)$ is the space of all measurable
bounded functions on $\X.$ For the second equality, note that in
the proof of Proposition \ref{res-21}, taking the test functions
$F(z,x)$ bounded, $z$-continuous and $x$-measurable (instead of
$x$-continuous), does not change anything except that in the
expression of the rate function  $\sup_{F}\{\langle F,q\rangle -
\psi(F)\},$ the sup is taken over this larger space instead of
$\CZX.$ As the rate function is unique, the sup over these two
spaces is the same. A similar argument in the present situation
leads to $\sup_{f\in \CX}=\sup_{f\in B(\X)}.$ Finally, choosing
$f\in B(\X)$ of the form $f=g\circ \beta$ with $g$ in $\CZ$ gives
us $
  \Psi^*(Q) \geq \sup_{g\in \CZ}\{\langle g,\beta\push Q\rangle-\langle g,\mu \rangle\}
$ and one concludes as in the previous proof.
\end{remark}

\section{Large deviations of a doubly indexed sequence of random measures. Preliminary results}
\label{sec-doubly indexed}

We keep the abstract Polish spaces $\Z$ and $\X$ of Section
\ref{sec-simply indexed}, as well as the triangular array
$(\zni\in\Z; 1\leq i\leq n, n\geq 1)$ which satisfies
(\ref{eq-mun}). For each $k\geq 1,$ we consider a Feller system
$(P^k_z\in\PX; z\in\Z)$ of probability laws on $\X$ and a
triangular array of \emph{independent} $\X$-valued random
variables $(\Xkni; 1\leq i\leq n, n\geq 1)$ where for each index
$(n,i)$ the law of $\Xkni$ is $P^k_\zni.$ This means that for all
$k,n\geq 1,$
\begin{equation*}
    \mathcal{L}aw(\Xkni; 1\leq i\leq n)=\otimes_{1\leq i\leq n}
    P^k_\zni.
\end{equation*}
The main result of the next Section \ref{sec-doubly main} states
the \knldp\ in $\PX$ for
\begin{equation*}
    \Lkn=\frac{1}{n} \sum_{i=1}^n \delta_\Xkni, \quad k,n\geq
    1.
\end{equation*}
As in Section \ref{sec-simply indexed}, this LDP will be obtained
by means of the contraction principle applied to some LDP for the
 $\PZX$-valued random variables
\begin{equation*}
    \Kkn=\frac{1}{n} \sum_{i=1}^n \delta_{(\zni,\Xkni)}, \quad k,n\geq
    1.
\end{equation*}
The main result of the present section is Theorem \ref{res-23}. It
states the \knldp\ for $\{\Kkn\}_{k,n\geq 1}.$

\par\medskip
We also assume that for each $z\in\Z,$ $(P^k_z)_{k\geq 1}$ obeys
some \kldp\ in $\X$ with rate function $J_z.$ This means that for
each $k$ and all measurable subset $A$ of $\X$
\begin{eqnarray*}
  -\inf_{x\in \inter A} J_z(x) &\leq& \lik\klog (A) \\
\nonumber  &\leq& \lsk\klog (A)\leq -\inf_{x\in \cl A} J_z(x)
\end{eqnarray*}
where $\inter A$ and $\cl A$ are the interior and the closure of
$A$ in $\X.$ This plays the part of Cram\'er's theorem and its
transformations at Section \ref{sec-main results}, see Examples
\ref{ex-02} with $J_z(x)=\cY(x_1-x_0),$ $x=(x_0,x_1)\in\Rdd,$
$z\in\Rd$ if $x_0=z$ and $+\infty$ otherwise.

\subsection{Preliminary results}
Before proving the \knldp\ for $\{\Kkn\}$ at Theorem \ref{res-23},
we need some preliminary results. The following lemma is Corollary
\ref{app-A2}, its detailed proof is given at Section
\ref{sec-appendix}.

\begin{lemma}\label{res-06}
 Let  $(\mathcal{F},\|\cdot\|)$ be a normed vector space and $\mathcal{Q}$
 be its dual space. Let $\lambda,$ $\lambda_k,$ $k\geq 1$ be
real-valued convex functions on $\mathcal{F}$ such that
\begin{itemize}
    \item[(a)] $\limk \lambda_k(F)=\lambda(F)$ for all $F\in
    \mathcal{F}$ and
    \item[(b)] there exists $c>0$ such that $\sup_{k\geq 1} |\lambda_k(F)|\leq
    c(1+\|F\|)$ for all $F\in \mathcal{F}.$
\end{itemize}
Then, the convex conjugates $\lambda_k^*$ of $\lambda_k,$
$\Gamma$-converge to the convex conjugate  $\lambda^*$ of
$\lambda:$
$$
\Glimk \lambda_k^*(q)=\lambda^*(q)
$$
for all $q\in\mathcal{Q},$ with respect to the $\ast$-weak
topology $\sigma(\mathcal{Q},\mathcal{F}).$
\end{lemma}

The following lemma is proved in \cite{Leo05a}.

\begin{lemma}\label{res-07}
 Suppose that for all $k\geq 1,$ $\{\mu^k_n\}_{n\geq 1}$ obeys a
 weak \nldp\ with rate function $kI^k$ and also
 suppose that the sequence $(I^k)_{k\geq 1}$ $\Gamma$-converges to
 some function $I.$ Then, $\{\mu^k_n\}_{k,n\geq 1}$ obeys a weak
 \knldp\ with rate function $I.$
\end{lemma}
\proof See \cite{Leo05a}.\endproof

We define for each $k$ and all $F\in \CZX,$
\begin{eqnarray}\label{eq-lambda}
 \lambda_k(F)&=&  \frac{1}{k} \IZ \log \langle
 e^{kF_z},P^k_z\rangle\,\mu(dz)\nonumber\\
        \lambda(F)
    &=& \IZ \supx \{F(z,x)-J_z(x)\}\,\mu(dz).
\end{eqnarray}
Note that $z\mapsto \supx \{F(z,x)-J_z(x)\}$ is measurable since
it is the pointwise limit of continuous functions: see
(\ref{eq-05}) below, so that $\lambda(F)$ is well-defined.
\\
Observe that $\lambda_k$ is a normalized version of the function
$\psi$ defined at (\ref{eq-psi}).

\begin{lemma}\label{res-34} We assume that
for each $z\in\Z,$ $(P^k_z)_{k\geq 1}$ obeys the \kldp\ in $\X$
with the \emph{good} rate function $J_z.$ Then, for all $F\in\CZX$
we have
\begin{align}
    \label{eq-05}
&\limk\frac{1}{k} \log \langle
 e^{kF_z},P^k_z\rangle=
    \supx \{F(z,x)-J_z(x)\},\\
    \label{eq-06}
&\limk \lambda_k(F)= \lambda(F)\quad \textrm{and}\\
\label{eq-08} &\sup_k |\lambda_k(F)|\leq \|F\|,\quad
|\lambda(F)|\leq \|F\|,\quad \forall F\in\CZX.
\end{align}

The functions $\lambda_k$ and $\lambda$ are convex and
$\sigma(\CZX,\CZX')$-\lsc..
\end{lemma}

 \proof
Thanks to the assumption on $(P^k_z)_{k\geq 1},$ by Varadhan's
integral lemma (see \cite{DZ}, Theorem 4.3.1), as $F_z$ is
continuous and bounded and $J_z$ is assumed to be a good rate
function, for all $z$ we have (\ref{eq-05}).

As for all $k\geq 1,$ $z\in\Z$ and $F\in \CZX,$ we have $
|\frac{1}{k}\log \langle e^{kF_z},P^k_z\rangle| \leq \|F\|, $ with
(\ref{eq-05}) we see that
\begin{equation}\label{eq-27}
   \left| \supx \{F(z,x)-J_z(x)\}\right|\leq \|F\|.
\end{equation}
These estimates allow us to apply Lebesgue dominated convergence
theorem to obtain (\ref{eq-06}) and (\ref{eq-08}).

 For each $k,$  $\lambda_k$ is convex since $f\mapsto
\log\langle e^{kf},P^k_z\rangle$ is convex as a log-Laplace
transform and $\mu$ is a nonnegative measure. As a pointwise limit
of convex functions, $\lambda$ is also convex.

The convex functions $\lambda^k$ and $\lambda$ are
$\sigma(\CZX,\CZX')$-\lsc\ if and only if they are
$\|\cdot\|$-\lsc\ on $\CZX.$ But, because of (\ref{eq-08}), these
convex functions are $\|\cdot\|$-continuous on the whole space
$\CZX.$ A fortiori, they are \lsc.
\endproof

\subsection{The \knldp\ for $\{\Kkn\}$}

Let us introduce the convex conjugate of $\lambda:$
\[
\lambda^*(q)=\sup_{F\in \CZX}\left\{\langle q,F\rangle -\IZ \supx
\{F(z,x)-J_z(x)\}\,\mu(dz)\right\}, \quad q\in \CZX'.
\]
It will appear during the proof of Theorem \ref{res-23} that it is
the rate function of the \knldp\ satisfied by $\{\Kkn\}_{k,n\geq
1}.$

\begin{thm}\label{res-23}
Suppose that
\begin{enumerate}
    \item $(\mu_n)_{n\geq 1}$ converges to $\mu$ in $\PZ,$
    \item for each $k\geq 1,$ $(P^k_z; z\in\Z)$ is a Feller
    system in the sense of Definition \ref{def-Feller},
    \item  for each $z\in\Z,$ $(P^k_z)_{k\geq 1}$ obeys the \kldp\
    in $\X$ with the good rate function $J_z.$
\end{enumerate}

Then $\{K^k_n\}_{k,n\geq 1}$ obeys the \knldp\ in $\PZX$ with the
 affine good rate function
\begin{equation}\label{eq-24}
i(q):=\left\{%
\begin{array}{ll}
    \IZX J_z(x)\,q(dzdx) & \hbox{if }q_\Z=\mu \\
    +\infty & \hbox{otherwise} \\
\end{array}%
\right.,\quad q\in\PZX.
\end{equation}
\end{thm}

 \proof
The framework of the proof is the same as Proposition
\ref{res-21}'s one, but it is technically more demanding.
\\
For all $k,n\geq 1$ and all $F\in \CZX,$ the normalized
log-Laplace transform of $K^k_n$ is defined by
\begin{equation*}
  \lambda_k^n(F)
  :=  \frac{1}{kn} \log \E \exp(kn\langle F,K^k_n\rangle)= \frac{1}{k} \IZ \log \langle
    e^{kF_z},P^k_z\rangle\,\mu_n(dz).
\end{equation*}
For fixed $k,$ considering the limit as $n$ tends to infinity and
taking assumptions (1) and (2) into account gives
\begin{equation*}
\limn \lambda_k^n(F) =  \lambda_k(F).
\end{equation*}
By Corollary \ref{res-201} we see that for all $k,$
$\{K^k_n\}_{n\geq 1}$ obeys the \nldp\ in $\CZX'$ with the rate
function $\lambda_k^*/k.$

\par\medskip\noindent
Because of Lemma \ref{res-34} and Lemma \ref{res-06} applied with
$\mathcal{F}=\CZX$ and $\mathcal{Q}=\CZX',$ the pointwise
convergence (\ref{eq-06}) and the estimate (\ref{eq-08}) imply
that
\begin{equation}\label{eq-31}
    \Glimk \lambda_k^*=\lambda^*
\end{equation}
in $\CZX'.$
\\
By Lemma \ref{res-07}, this $\Gamma$-convergence
 implies that $\{K^k_n\}_{k,n\geq
1}$ obeys a \emph{weak} \knldp\ in $\CZX'$ with the rate function
$\lambda^*.$ It is proved at Lemma \ref{res-24} below that
\begin{equation}\label{eq-217}
    \{\lambda^*<+\infty\}\subset\PZX.
\end{equation}
A fortiori, $\{\lambda^*<+\infty\}$ is included in the strong unit
ball
$$
\UZX =\left\{q\in \CZX'; \|q\|^*:=\sup_{F\in \CZX,
\|F\|\leq 1}\langle q,F\rangle\leq 1\right\}
$$
of $\CZX'$ which is $\sigma(\CZX',\CZX)$-compact (Banach-Alaoglu
theorem). Consequently, $\{K^k_n\}_{k,n\geq 1}$ obeys a
\emph{strong} \knldp\ in $\UZX $ with the topology $\sigma(\UZX
,\CZX)$ and the rate function $\lambda^*.$ With (\ref{eq-217})
again, we obtain that $\{K^k_n\}_{k,n\geq 1}$ obeys the \knldp\ in
$\PZX$ with the rate function $\lambda^*.$

\par\medskip\noindent
 Let us show that the restriction of $\lambda^*$ to $\PZX$ has
$\sigma(\PZX,\CZX)$-compact level sets.
 As a convex conjugate, $\lambda^*$ is $\sigma(\CZX',\CZX)$-\lsc.
 Therefore, for all real $a,$ $\{\lambda^*\leq a\}$ is
 $\sigma(\CZX',\CZX)$-closed. But, (\ref{eq-217}) implies that $\{\lambda^*\leq a\}$
 is included in the $\sigma(\CZX',\CZX)$-compact unit ball $\UZX .$ Hence,
 $\{\lambda^*\leq a\}$ is $\sigma(\CZX',\CZX)$-compact and by (\ref{eq-217}) again,
  the restriction of $\lambda^*$ to $\PZX$ is $\sigma(\PZX,\CZX)$-inf-compact.

\par\medskip\noindent
 Finally, it will be proved at
Proposition \ref{res-31} that the restriction of $\lambda^*$ to
$\PZX$ is $i.$ This completes the proof of the theorem.
\endproof

\subsection{Identification of the rate function $\lambda^*$}
It remains to show that $\lambda^*=i.$ This is the most technical
part of the paper.

\begin{lemma}\label{res-24}
Under the assumptions of Lemma \ref{res-34}, the following
statements hold true.
\begin{itemize}
    \item[(a)]
    For all $q\in \CZX',$ $\lambda^*(q)<\infty$
    implies that $q\in\PZX.$
    \item[(b)]
    For all $q\in\PZX,$ $\lambda^*(q)<\infty$
    implies that $q_\Z=\mu.$
\end{itemize}
\end{lemma}

 \proof
It is similar to the proof of Lemma \ref{res-202}. As in Lemma
\ref{res-202},  the $\|\cdot\|$-continuity of $q$ doesn't play any
role, see Remark \ref{rem-01}. Let $q\in \CZX'$ be such that
\begin{equation*}
    \sup_{F\in \CZX}\{\langle F,q\rangle-\lambda(F)\}=\lambda^*(q)<\infty.
\end{equation*}

An inspection of Lemma \ref{res-202}'s proof shows that, to prove
that $q\in\PZX,$ it is enough to check that $\lambda$ satisfies
\begin{itemize}
    \item[(i)] for all $a\leq 0$ and all nonnegative $F_o\in \CZX,$ $\lambda( aF_o)\leq 0$
    \item[(ii)] for any constant function $F\equiv c\in \mathbb{R},$ $\lambda(c\1)=c$
    \item [(iii)] for any sequence $(F_n)_{n\geq 1}$ in $\CZX$ such that $F_n\geq 0$
    for all $n$ and $(F_n(z,x))_{n\geq 1}$ decreases to zero for each $(z,x)\in\ZX,$ we have,
    $\lim_{n\rightarrow\infty}\lambda(aF_n)=0, $ for all $a\geq 0.$
\end{itemize}

 \noindent(i)\quad
 As $J_z(x)\geq 0$ for all $z$ and $x,$ and
$\mu\geq 0,$ we have $\lambda( aF_o)\leq 0$ for all $a\leq 0$ and
all nonnegative $F_o\in \CZX.$

 \par\medskip\noindent (ii)\quad
 As $\infx J_z(x)=0$ for all $z$ and $\mu$ is a probability measure, for any
constant function $F\equiv c\in \mathbb{R},$ we have
$\lambda(c\1)=c.$

 \par\medskip\noindent (iii)\quad
By Lemma \ref{res-26} below, for all $z\in\Z,$ $\big(\supx
\{F_n(z,x)-J_z(x)\}\big)_{n\geq 1}$ is a decreasing sequence and
$\lim_{n\rightarrow\infty}\supx \{F_n(z,x)-J_z(x)\}=0.$ As
$\left|\supx \{F_n(z,x)-J_z(x)\}\right|\leq
\sup_{z,x}|F_1(z,x)|<\infty$ for all $n$ and $z,$ one can apply
the dominated convergence theorem to obtain that $
\lim_{n\rightarrow\infty}\lambda(aF_n)=0, $ for all $a\geq 0.$

This completes the proof of statement (a).
\par\medskip\noindent
Let us prove (b). Choosing $F(z,x)=g(z)$ not depending on $x$ with
$g\in\CZ$ in the expression of $\lambda^*,$ one sees that for all
$q\in\PZX$
\begin{eqnarray*}
  \lambda^*(q) &\geq& \sup_{g\in\CZ}\{\langle g, q_\Z\rangle-\langle g,\mu \rangle\} \\
   &=& \left\{%
\begin{array}{ll}
    0, & \hbox{if }q_\Z=\mu \\
    +\infty, & \hbox{otherwise} \\
\end{array}%
\right.
\end{eqnarray*}
which gives the announced result and completes the proof of Lemma
\ref{res-24}.
\endproof

The very technical result of this section is the following
Proposition \ref{res-31}. During its proof, we need some lemmas
whose statements are included in the body of the proof. The proofs
of these lemmas are postponed to the next subsection
\ref{sec-lemmas}.

\begin{prop}\label{res-31}
For all $q\in\PZX,$ $\lambda^*(q)=i(q).$
\end{prop}
\proof
 Thanks to Lemma \ref{res-24}-b, to prove that
$ \lambda^*=i, $ we have to show that $\lambda^*(q)=\IZX
J_z(x)\,q(dzdx)$ for all $q\in\PZX$ such that $q_\Z=\mu$ or
equivalently such that
\begin{equation}\label{eq-25}
q(dzdx)=\mu(dz)q^z(dx)
\end{equation}
where
$$
q^z(dx)=q(X\in dx\mid Z=z).
$$
For such a $q$ we have
\begin{eqnarray*}
  \lambda^*(q)
  &=& \sup_{F\in\CZX}\IZ [\langle F_z,q^z\rangle-\supx \{F_z(x)-J_z(x)\}]\,\mu(dz) \\
  &\leq& \IZ \sup_{f\in\CX}[\langle f,q^z\rangle-\supx \{f(x)-J_z(x)\}]\,\mu(dz) \\
    &\stackrel{(a)}{=}& \IZX J_z(x)\, q^z(dx)\mu(dz)\\
    &\stackrel{(b)}{=}& i(q)
\end{eqnarray*}
where equality (a) is given at the following Lemma \ref{res-25a}
and equality (b) follows from (\ref{eq-25}).

\begin{lemma}\label{res-25a} Let $J$ be a $[0,+\infty]$-valued \lsc\
function on $\X.$
    For  all $Q\in\PX,$ we have
    $$
\sup_{f\in\CX}\left\{\IX f\,dQ-\supx (f(x)-J(x))\right\}=\IX
J\,dQ.
$$
\end{lemma}
The proof of this lemma is put back after the proof of the present
proposition.\par\medskip

Note that $z\mapsto \langle J_z,q^z\rangle$ is measurable since
$z\mapsto J_z(x)$ is assumed to be continuous for all $x$ and
$z\mapsto q^z$ is a regular version of the desintegration of $q.$

It remains to show the converse inequality: $\lambda^*(q)\geq
i(q)$ for all $q$ satisfying (\ref{eq-25}). As a first step, we
would like to invert a sup and an integral to obtain
\begin{eqnarray}\label{eq-310}
\nonumber \lambda^*(q)  &=& \sup_{F\in \CZX)}\IZ [\langle F_z,q^z\rangle-\supx \{F_z(x)-J_z(x)\}]\,\mu(dz) \\
  &=& \IZ \sup_{f\in\CX}[\langle
f,q^z\rangle-\supx \{f(x)-J_z(x)\}]\,\mu(dz)
\end{eqnarray}
As a first step, we are going to prove this equality under the
restrictive assumption that $\X$ is compact. Its proof relies on
the following result which is due to R.~T.~Rockafellar (see
\cite{Roc71}, Theorem 2).

\begin{lemma}\label{res-32}
Let $(\Z,\mu)$ be a measure space such that $\mu$ is
$\sigma$-finite. Let $L$ be a \emph{decomposable space} (see below
for the definition) of measurable functions $F$ on $\Z$ with their
values in a Polish space $\mathcal{Y}$ equipped with its Borel
$\sigma$-field. Let $\theta: \Z\times \mathcal{Y}\to
[-\infty,\infty)$ be such that
\begin{itemize}
    \item[-] $\theta$ is jointly measurable
    \item[-] $\theta$ is not identically equal to $-\infty$ and
    \item[-] $y\mapsto\theta(z,y)$ is upper semicontinuous for all
    $z\in\Z.$
\end{itemize}
In this case, one says that $-\theta$ is \emph{normal}. Suppose in
addition that there exist some $F_1\in L$ and some $u_1\in
L^1(\mu)$ such that $\theta(z,F_1(z))\geq u_1(z)$ for $\mu$-almost
every $z$ in $\Z.$ Then, $z\mapsto \sup_{y\in \mathcal{Y}}
\theta(z,y)$ is measurable and
\[
\sup_{F\in L} \IZ \theta(z,F(z))\,\mu(dz)
 =\IZ \sup_{y\in \mathcal{Y}} \theta(z,y)\, \mu(dz)\in
 (-\infty,\infty].
\]
\end{lemma}

\begin{definition}
The space $L$ is said to be \emph{decomposable} if, whenever $F$
belongs to $L$ and $F_o:\Z_o\to\mathcal{Y}$ is a bounded
measurable function on a measurable set $\Z_o\subset\Z$ of finite
measure, the function $z\mapsto\1_{z\in\Z_o}
F_o(z)+\1_{z\not\in\Z_o}F(z)$ also belongs to $L.$
\end{definition}

In order to obtain (\ref{eq-310}), we would like to apply this
lemma with
\begin{itemize}
    \item $\mathcal{Y}=\CX$ equipped with the topology of uniform convergence,
    \item $\theta(z,f)=\langle
    f,q^z\rangle-\supx \{f(x)-J_z(x)\}$  and
    \item $L=C_b(\Z,\CX)\simeq\CZX.$
\end{itemize}

Unfortunately, two troubles occur.
\\
\textsl{Trouble 1}: If $\X$ is not compact, $\mathcal{Y}=\CX$ is
\emph{not} separable and fails to be a Polish space as required in
the lemma. On the other hand, if $\X$ is compact, $\CX$ is Polish.
\\
\textsl{Trouble 2}: The space $\CZX\simeq C_b(\Z,\CX)$ is
\emph{not} decomposable. On the other hand, the space $B(\Z,\CX)$
of all bounded and measurable functions $F: z\in\Z\mapsto
F_z\in\CX$ is decomposable.
\\
Note that when $\X$ is compact, as $\CX$ is separable, we have $
B(\Z,\CX)\simeq \mathcal{F} $ where $\mathcal{F}$ is the space of
all the functions on $\ZX$ which are bounded, $x$-continuous and
$z$-measurable; such functions are jointly measurable.
\\
We are going to apply Lemma \ref{res-32} with
\begin{itemize}
    \item $\mathcal{Y}=\CX$  and $\X$ a compact Polish set,
    \item $\theta(z,f)=\langle
    f,q^z\rangle-\supx \{f(x)-J_z(x)\}$ for all $z\in\Z$ and $f\in\CX,$ where $q\in\PZX$ is
    fixed and satisfies (\ref{eq-25}) and
    \item $L=\mathcal{F}.$
\end{itemize}

As $f\mapsto \theta(z,f)$ is continuous for all $z$ and $z\mapsto
\theta(z,f)$ is measurable for all $f,$ $\theta$ is jointly
measurable. Taking $f=0$ gives $\theta(z,0)=0>-\infty$ for all
$z,$ so that $\theta$ shares all the normality conditions of the
lemma.
\\
Choosing the functions $F_1=0\in L$ and $u_1=0\in L^1(\mu)$ leads
us to $0=\theta(z,F_1(z))\geq u_1(z)=0$ for every $z$ in $\Z.$

Therefore, we have shown that all the assumptions of Lemma
\ref{res-32} are met so that
\begin{eqnarray}\label{eq-26}
\nonumber  && \sup_{F\in \mathcal{F}}\IZ [\langle F_z,q^z\rangle-\supx \{F_z(x)-J_z(x)\}]\,\mu(dz) \\
           &=& \IZ \sup_{f\in\CX}[\langle f,q^z\rangle-\supx \{f(x)-J_z(x)\}]\,\mu(dz)
\end{eqnarray}
whenever $\X$ is a compact Polish space.

To obtain (\ref{eq-310}), it remains to prove that for all $q$
with $q_\Z=\mu,$
\begin{equation}\label{eq-32}
 \lambda^*(q)
= \sup_{F\in \mathcal{F}}\IZ [\langle F_z,q^z\rangle-\supx
\{F_z(x)-J_z(x)\}]\,\mu(dz)
\end{equation}
Let us prove it without assuming that $\X$ is compact. Rather than
invoking an abstract approximation argument, we present a specific
proof of (\ref{eq-32}). Rewriting the above proof of Theorem
\ref{res-23} with $\CZX$ replaced with the space $B(\ZX)$ of
bounded measurable functions on $\ZX$ one gets the following
result.

\par\medskip
\textbf{A variant of Theorem \ref{res-23}.}\quad \textit{Assuming
(2) and (3) of Theorem \ref{res-23}, if Assumption (1) is
strengtnened by ``$(\mu_n)_{n\geq 1}$ converges to $\mu$ in $\PZ$
for the stronger topology $\sigma(\PZ,B(\Z))$", then
$\{K^k_n\}_{k,n\geq 1}$ obeys the \knldp\ in $\PZX$ with the
topology $\sigma(\PZX,B(\ZX))$ and the rate function
$\tilde{i}(q)= \sup_{F\in B(\ZX)}\IZ [\langle F_z,q^z\rangle-\supx
\{F_z(x)-J_z(x)\}]\,\mu(dz),$ if $q\in\PZX$ satisfies $q_\Z=\mu$
and $\tilde{i}(q)=+\infty$ otherwise.}

\par\medskip
For any $\mu\in\PZ,$ there exists a sequence of empirical measures
$(\mu_n)_{n\geq 1}$ as in (\ref{eq-mun}) which converges to $\mu$
with respect to the topology $\sigma(\PZ,B(\Z)).$ This can be seen
as a consequence of the almost sure convergence, as $n$ tends to
infinity, of the empirical measures $\frac{1}{n}\sum_{1\leq i\leq
n}\delta_{Z_i}$ of the $\mu$-iid sequence of $\Z$-valued random
variables $(Z_i)_{i\geq 1}$ towards $\mu$ for the topology
$\sigma(\PZ,B(\Z))$ which in turns is a corollary of the
strenghened version of Sanov's theorem  with the topology
$\sigma(\PZ,B(\Z))$ on a Polish space $\Z.$ With such a sequence
$(\mu_n)_{n\geq 1},$ by Theorem \ref{res-23} and its variant,
$\{K^k_n\}_{k,n\geq 1}$ obeys the \knldp\ in $\PZX$ with the rate
functions $\lambda^*$ and $\tilde{i}.$ As the rate function of a
LDP is unique in a regular space (for the double index version of
this known result, see \cite{Leo05a}), we have
$\lambda^*=\tilde{i}.$ It follows that for all $q$ with
$q_\Z=\mu,$
\begin{eqnarray*}
\lambda^*(q)&:=& \sup_{F\in \CZX}\IZ [\langle F_z,q^z\rangle-\supx
 \{F_z(x)-J_z(x)\}]\,\mu(dz)\nonumber\\
 &=& \sup_{F\in B(\ZX)}\IZ [\langle F_z,q^z\rangle-\supx \{F_z(x)-J_z(x)\}]\,\mu(dz)
\end{eqnarray*}
which implies the desired equality (\ref{eq-32}).

Thanks to (\ref{eq-26}) and (\ref{eq-32}), we have proved
(\ref{eq-310}) whenever $\X$ is compact. Nevertheless, the
identity (\ref{eq-310}) will not be used directly. We shall only
use (\ref{eq-32}) and a variant of (\ref{eq-26}).

\par\medskip
Now, we have to tackle the problem of relaxing the requirement
that $\X$ is compact. Let us take advantage of the tightness of
$q$ (it is a probability on a Polish space). This means that there
exists an increasing sequence $(\mathcal{K}_n^q)_{n\geq 1}$ of
compact subsets of $\ZX$ such that $q(\mathcal{K}_n^q)\geq 1-1/n$
for all $n\geq 1.$ As a continuous image of a compact set,
$\X_n^q:=\{x\in\X; (z,x)\in \mathcal{K}_n^q \textrm{ for some
}z\in\Z \}$ is a compact set. We also have $q(\Z\times\X_n^q)\geq
1-1/n$ for all $n.$ It follows that for $q_\Z$-almost every
$z\in\Z,$ $q^z$ is determined by the values $\langle f,q^z\rangle$
where $f$ describes the set $\bigcup_{n\geq 1}
\mathcal{C}(\X_n^q)$ where, for any measurable set $\X_o$ in $\X,$
we denote
\begin{equation}\label{eq-330}
    \mathcal{C}(\X_o)=\1_{\X_o}\CX=\{f: \X\to\mathbb{R};
f=\1_{\X_o}\tilde{f}, \textrm{ for some }\tilde{f}\in \CX\}.
\end{equation}

To see this, remark that for all measurable set $A$ in $\X$ such
that $A\cap (\cup_n\X_n^q)=\emptyset,$ we have $\IZ
q^z(A)\,\mu(dz)=q(\Z\times A)=\lim_{n\rightarrow\infty}q(\Z\times
(A\cap\X_n^q))=0.$

\par\medskip
We can now proceed with the proof of $\lambda^*(q)\geq i(q)$ for
all $q$ satisfying (\ref{eq-25}). For all such $q$ we have,

\begin{eqnarray*}
\lambda^*(q)
  &=& \sup_{F\in \mathcal{F}}\IZ [\langle F_z,q^z\rangle-\supx \{F_z(x)-J_z(x)\}]\,\mu(dz)\\
    &\stackrel{(a)}{\geq}& \sup_{n\geq 1} \sup_{F\in B\left(\Z,\mathcal{C}(\X_n^q)\right)}\IZ [\langle
                F_z,q^z\rangle-\supx \{F_z(x)-J_z(x)\}]\,\mu(dz)\\
  &\stackrel{(b)}{=}&  \sup_{n\geq 1} \IZ \sup_{f\in \mathcal{C}(\X_n^q)}[\langle
              f,q^z\rangle-\supx \{f(x)-J_z(x)\}]\,\mu(dz)\\
  &\stackrel{(c)}{=}&  \IZ \sup_{f\in \bigcup_{n\geq 1}\mathcal{C}(\X_n^q)}[\langle
              f,q^z\rangle-\supx \{f(x)-J_z(x)\}]\,\mu(dz)\\
 &\stackrel{(d)}{=}& \IZX J_z(x)\, q^z(dx)\mu(dx)\\
    &=& i(q)
\end{eqnarray*}
where the first equality is (\ref{eq-32}). The remaining series of
inequality and equalities needs to be justified. This will require
two more lemmas the proofs of which are postponed after the proof
of the present proposition.
\par\medskip

 $\bullet$ Inequality (a).
It is enough to show that for any function $F\in
B(\Z,\mathcal{C}(\X_o))$ with $\X_o$ a compact subset of $\X,$
there exists a sequence $(F^n)_{n\geq 1}$ in $\mathcal{F}$ such
that
\begin{eqnarray}\label{eq-214}
   && \IZ [\langle F_z,q^z\rangle-\supx \{F_z(x)-J_z(x)\}]\,\mu(dz)\\
 \nonumber  &=&\limn \IZ [\langle F_z^n,q^z\rangle-\supx \{F_z^n(x)-J_z(x)\}]\,\mu(dz).
\end{eqnarray}

Let us show that
$$
F^n_z(x):= \sup\{F_z(y)-nd(x,y); y\in\X\}
$$
does this job. For each $z,$ $(-F^n_z)_{n\geq 1}$ is the
Moreau-Yosida approximation of $-F_z,$ and it is a well-known
result (see \cite{Braides}, Section 1.7.3 for instance) that
\begin{itemize}
    \item[-] for all $z\in\Z,$ $x\mapsto F^n_z(x)$ is $n$-Lipschitz,
\end{itemize}
and for all $(z,x)\in\ZX,$
\begin{itemize}
    \item[-] $-\|F\|\leq F_z(x)\leq F^n_z(x)\leq \|F_z\|\leq \|F\|,$ where $\|\cdot\|$ stands for the uniform
    norm,
    \item[-] $(F^n_z(x))_{n\geq 1}$ is a decreasing sequence and
    \item[-] $\limn F^n_z(x)=F_z(x)$
\end{itemize}
For the last statement, note that it is necessary that $F_z$ is
upper semicontinuous on $\X.$ But, this is insured by the
assumption that $\X_o$ is closed and $F_z\in \mathcal{C}(\X_o).$

Now let us make sure that for any $x_o\in\X,$ $z\mapsto
F^n_z(x_o)$ is measurable. For all real $a,$ we have
\begin{eqnarray*}
  F_z^n(x_o)\leq a
  &\Leftrightarrow& \forall y\in\X, F_z(y)-nd(x_o,y)\leq a \\
   &\Leftrightarrow& \forall k\geq 1, F_z(x_k)-nd(x_o,x_k)\leq a
\end{eqnarray*}
where $\{x_k; k\geq 1\}$ is a countable dense subset of $\X$
(recall that $\X$ is Polish). This holds, since $y\mapsto
F_z(y)-nd(x_o,y)$ is continuous. It follows that $\{z\in\Z;
F_z^n(x_o)\leq a\}=\cap_{k\geq 1}\{z\in\Z;
F_z(x_k)-nd(x_o,x_k)\leq a\}.$ As $z\mapsto F_z(x_k)$ is
measurable for all $k,$ this proves the measurability of $z\mapsto
F^n_z(x_o).$ Therefore, $F^n$ belongs to $\mathcal{F}$ for all
$n\geq 1.$

With the estimate $-\|F\|\leq F_z(x)\leq F^n_z(x)\leq \|F_z\|\leq
\|F\|$ and the limit  $\limn F^n_z(x)=F_z(x),$ one can apply the
dominated convergence theorem to obtain that
\begin{equation}\label{eq-215a}
    \limn \IZ \langle F^n_z,q^z\rangle \,\mu(dz)=\limn \IZX
F^n\,dq=\IZX F\,dq=\IZ \langle F_z,q^z\rangle\,\mu(dz).
\end{equation}

Similarly, the limit
\begin{equation}\label{eq-215b}
\limn \IZ \supx \{F_z^n(x)-J_z(x)\}\,\mu(dz) = \IZ \supx
\{F_z(x)-J_z(x)\}\,\mu(dz)
\end{equation}
follows from the estimate (\ref{eq-27}) and the following lemma.

\begin{lemma}\label{res-26}
Let $J$ be an inf-compact $[0,\infty]$-valued function on $\X$ and
$(f_n)_{n\geq 1}$ a decreasing sequence of continuous bounded
functions on $\X$ which converges pointwise to some bounded upper
semicontinuous function $f.$ Then,
$\left(\supx\{f_n(x)-J(x)\}\right)_{n\geq 1}$ is a decreasing
sequence and
$$
\limn \supx\{f_n(x)-J(x)\}=\supx\{f(x)-J(x)\}.
$$
\end{lemma}
The proof of this lemma is put back after the proof of the present
proposition.\par\medskip

Finally, (\ref{eq-214}) follows from (\ref{eq-215a}) and
(\ref{eq-215b}).

\par\medskip
$\bullet$ Equality (b) is a variant of (\ref{eq-26}) applied with
the compact set $\X_n^q.$

\par\medskip
$\bullet$ Equality (c). If the sequence
$\mathcal{C}(\X_n^q)_{n\geq 1}$ were increasing, equality (c)
would be a direct consequence of the monotone convergence theorem.
Nevertheless, this is almost the case since, for any pair of
closed subsets $\X_o$ and $\X_1$ of $\X$ such that $\X_o\subset
\X_1,$ any function $f\in \mathcal{C}(\X_o)$ can be approximated
pointwise by a uniformly bounded decreasing sequence $(f_n)$ in
$\mathcal{C}(\X_1)$ such that $\limn \supx \{f_n(x)-J_z(x)\}=\supx
\{f(x)-J_z(x)\}.$ One proves this, exactly as for inequality (a),
by means of a Moreau-Yosida approximation and Lemma \ref{res-26}.
With this in hand, equality (c) follows from the monotone
convergence theorem.

\par\medskip
$\bullet$ Equality (d). This equality is a consequence of the
following lemma.

\begin{lemma}\label{res-25b}
Let $J$ be a $[0,+\infty]$-valued \lsc\ function on $\X.$

If $\CX$ in Lemma \ref{res-25a} is replaced with the set $
\mathcal{G}_Q=\bigcup_{n\geq 1} \mathcal{C}(\X_n^Q) $ where
$(\X_n^Q)_{n\geq 1}$ is an increasing sequence of closed subsets
of $\X$ such that $\limn Q(\X_n^Q)=1,$ then we still have
$$
\sup_{f\in\mathcal{G}_Q}\left\{\IX f\,dQ-\supx
(f(x)-J(x))\right\}=\IX J\,dQ.
$$
\end{lemma}
The proof of this lemma is put back after the proof of the present
proposition.\par\medskip

Note that we have already remarked that for $q_\Z$-almost every
$z\in\Z,$ $q^z$ is determined by the values $\langle f,q^z\rangle$
where $f$ describes the set $\bigcup_{n\geq 1}
\mathcal{C}(\X_n^q).$ One obtains equality (d) by means of Lemma
\ref{res-25b}, with $\mathcal{G}_{q^z}=\bigcup_{n\geq
1}\mathcal{C}(\X_n^q),$ for all $z\in\Z.$

\par\medskip
We have proved that $\lambda^*\geq i$ and this completes the proof
of the proposition.
\endproof

\par\medskip
\textit{A comment on this proof.}\ One could think of replacing
the spaces $\mathcal{C}(\X_n^q)$ defined by (\ref{eq-330}) with
the smaller spaces $\hat{\mathcal{C}}(\X_n^q)$ of continuous
functions on $\X$ with their support in $\X_n^q.$ This clearly
provides an increasing sequence and simplifies the proof of
equality (c). But unfortunately, equality (a) doesn't work anymore
since $\hat{\mathcal{C}}(\X_o)$ reduces to the null space when the
compact set $\X_o$ has an empty interior (a common feature in
infinite dimension).

\subsection{Proofs of the lemmas}\label{sec-lemmas}
We go on with the proofs of Lemmas \ref{res-25a} , \ref{res-26}
and \ref{res-25b}.

\par\medskip
 \proof[Proof of Lemmas \ref{res-25a} and \ref{res-25b}]
Lemma \ref{res-25a} is a particular case of Lemma \ref{res-25b},
we only prove Lemma \ref{res-25b}.

As,
    $\sup_{R\in \PX}\langle f-J,R\rangle
    \leq \supx \{f(x)-J(x)\}
    =\supx\langle f-J,\delta_{x}\rangle
    \leq\sup_{R\in \PX}\langle f-J,R\rangle,$ we have
     $\supx \{f(x)-J(x)\}=\sup_{R\in \PX}\langle f-J,R\rangle.$
     Therefore,
\begin{eqnarray*}
  \sup_{f\in\mathcal{G}_Q}\left\{\langle f,Q\rangle-\sup_{x\in\X}\{f(x)-J(x)\}\right\}
  &=& \sup_{f\in\mathcal{G}_Q}\left\{\langle f,Q\rangle-\sup_{R\in \PX}\langle f-J,R\rangle\right\} \\
  &=& \langle J,Q\rangle +\sup_{f\in\mathcal{G}_Q}\left\{\langle f-J,Q\rangle-\sup_{R\in \PX}\langle f-J,R\rangle\right\} \\
  &\leq& \langle J,Q\rangle
\end{eqnarray*}
where the last inequality holds since $Q\in\PX.$

\par\medskip
Now, let's prove the converse inequality. As $J$ is a \lsc\
function which is bounded below, it is the pointwise limit of an
increasing sequence $(\tilde{J}_n)_{n\geq 1}$ in $\CX$: once
again, the Moreau-Yosida approximation:
$\tilde{J}_n(x)=\inf\{J(y)+nd(x,y); y\in\X\}.$

Let us define $J_n(x)=\1_{\X_n^Q}(x)(0\vee\tilde{J}_n(x)\wedge n)$
for all $x$ and $n.$ As $(\X_n^Q)_{n\geq 1}$ is an increasing
sequence of sets, $(J_n)_{n\geq 1}$ is an increasing sequence of
functions such that for all $n,$ $J_n$ is in
$\mathcal{C}(\X_n^Q).$ We have
\begin{eqnarray*}
   \sup_{f\in\mathcal{G}_Q}\left\{\langle f,Q\rangle-\sup_{x\in\X}\{f(x)-J(x)\}\right\}
   &\stackrel{(a)}{\geq}& \sup_{n\geq 1} \left(\langle J_n,Q\rangle -\sup_{x\in\X}\{J_n(x)-J(x)\}\right) \\
   &\stackrel{(b)}{\geq}& \sup_{n\geq 1}\IX J_n\,dQ \\
    &\stackrel{(c)}{\geq}& \sup_{k\geq 1}\sup_{n\geq k}\int_{\X_k^Q} (0\vee \tilde{J}_n\wedge n)\,dQ, \\
      &\stackrel{(d)}{=}& \sup_{k\geq 1}\int_{\X_k^Q} J\,dQ, \\
   &\stackrel{(e)}{=}& \IX J\,dQ
\end{eqnarray*}
where inequality (a) holds since $J_n\in\mathcal{C}(\X_n^Q),$
inequality (b) follows from $J_n\leq J,$ equality (c) holds since
the sequence $(\X_n^Q)$ is increasing, equality (d) follows from
the monotone convergence theorem and equality (e) follows from the
monotone convergence theorem together with
$\lim_{k\rightarrow\infty}Q(\X\setminus \X_k^Q)=0.$ This completes
the proof of the lemmas.
 \endproof

 \proof[Proof of Lemma \ref{res-26}]
Changing sign and denoting $g_n(x)= J(x)-f_n(x),$
$g(x)=J(x)-f(x),$ we want to prove that $\limn\infx g_n(x)=\infx
g(x).$

We see that $(g_n)_{n\geq 1}$ is an increasing sequence of \lsc\
functions. It follows by the Proposition 5.4 of \cite{DalMaso}
that it is a $\Gamma$-convergent sequence and
\begin{equation}\label{eq-216a}
    \Glimn g_n=\limn g_n=g.
\end{equation}
Let us admit for a while that there exists some compact set $K$
which satisfies
\begin{equation}\label{eq-216b}
     \infx g_n(x)= \inf_{x\in K}g_n(x)
\end{equation}
for all $n.$ This and the convergence (\ref{eq-216a}) allows to
apply Theorem 7.4 of \cite{DalMaso} to obtain $\limn \infx g_n(x)=
\infx \Glimn g_n(x)=\infx g(x)$ which is the desired result.

\par\medskip
It remains to check that (\ref{eq-216b}) is true. Let $x_*\in\X$
be such that $J(x_*)<\infty$ (if $J\equiv +\infty,$ there is
nothing to prove). Then, $\infx g_n(x)\leq
g_n(x_*)=J(x_*)-f_n(x_*)\leq J(x_*)-f(x_*) \leq J(x_*)-\infx
f(x).$ On the other hand, for all $x$ and $n,$ $f_n(x)\leq f_1(x)
\leq A:= \sup f_1.$ Let $B:= A+1+J(x_*)-\infx f(x).$ For all $x$
such that $J(x)>B,$ we have $g_n(x)>B-\supx f_n(x)\geq B-A\geq
J(x_*)-\infx f(x)+1.$ We have just seen that for all $n,$
\begin{eqnarray*}
    \infx g_n(x)&\leq& J(x_*)-\infx f(x) \\
    \inf_{x; J(x)>B} g_n(x)&\geq& J(x_*)-\infx f(x)+1
\end{eqnarray*}
This proves (\ref{eq-216b}) with the compact level set $K=\{J\leq
B\}$ and completes the proof of the lemma.
\endproof

\section{Large deviations of a doubly indexed sequence of random measures. Main results}
\label{sec-doubly main}

 Theorem \ref{res-23} states a
\knldp\ for $K_n^k=\frac{1}{n} \sum_{i=1}^n \delta_{(\zni,\Xkni)}$
but we are mostly interested in the \knld\  in $\PX$ of $
    L_n^k=\frac{1}{n} \sum_{i=1}^n \delta_\Xkni.
$ It will easily follow from Theorem \ref{res-23} and the
contraction principle. Let us denote
\[
P^k(dx)= \IZ P_z^k(dx)\,\mu(dz)\in\PX, k\geq 1.
\]
\begin{thm}\label{res-33}
Suppose that
\begin{enumerate}
    \item $(\mu_n)_{n\geq 1}$ converges to $\mu$ in $\PZ,$
    \item for each $k\geq 1,$ $(P^k_z; z\in\Z)$ is a Feller
    system in the sense of Definition \ref{def-Feller},
    \item  for each $z\in\Z,$ $(P^k_z)_{k\geq 1}$ obeys the \kldp\
    in $\X$ with the good rate function $J_z.$
\end{enumerate}
Then the following statements hold true.
\begin{enumerate}
 \item[(a)]
    $\{L_n^k\}_{k,n\geq 1}$ obeys the \knldp\ in $\PX$ with the good rate
    function $I$ which is defined for all $Q\in\PX$ by
    \begin{equation}\label{eq-35}
    I(Q)= \inf\left\{\IZX  J_z(x) \,\mu(dz)\Pi_z(dx); (\Pi_z)_{z\in\Z} : \IZ \Pi_z\,\mu(dz)=Q\right\}
    \end{equation}
    where the transition kernels $z\in\Z\mapsto \Pi_z\in\PX$ are measurable.
 \item[(b)]
    Another representation of this rate function is
    \begin{equation*}
    I(Q)=\sup_{f\in\CX}\left\{\IZ \mathcal{S}f(z)\,\mu(dz)-\IX f(x)\,Q(dx)\right\}, Q\in\PX
    \end{equation*}
    where $\mathcal{S}f(z)$ is defined for all $z\in\Z$ by
    \begin{equation*}
    \mathcal{S}f(z)=\infx\{J_z(x)+f(x)\}.
\end{equation*}
 \item[(c)]
    If $I(Q)<+\infty,$ there exists a (possibly not unique) kernel $(\Pi^*_z)_{z\in\Z}$
    which realizes the infimum in (\ref{eq-35}).
 \item[(d)]
    If for each $k$ the Feller system
$(P_z^k)_{z\in\Z}$ satisfies
\begin{equation}\label{eq-301}
     P_z^k=P^k (\cdot\mid \beta(X)=z)
\end{equation}
    for $\mu$-almost every $z\in\Z$ and some continuous function
    $\beta: \X\to\Z,$ we have
    \begin{equation}\label{eq-41}
I(Q)=\left\{\begin{array}{ll}
     \IX J_{\beta(x)}(x)\,Q(dx)
      & \textrm{if }\beta\push Q=\mu\\
      +\infty & \textrm{otherwise} \\
    \end{array}\right.,\quad Q\in\PX.
\end{equation}
\end{enumerate}
\end{thm}

\noindent The dual space $\CX'$ of $(\CX,\|\cdot\|)$ is equipped
with the $\ast$-weak topology $\sigma(\CX',\CX),$ see Section
\ref{sec-conventions}.

\proof
 Let us prove (a). As $L_n^k$ is the $\X$-marginal of $K_n^k$ and
$\{K_n^k\}$ obeys the \knldp\ with the good rate function $i,$ the
statement (a) follows from an obvious extension to the double
index setting of the contraction principle (see \cite{Leo05a}):
$\{L_n^k\}$ obeys the \knldp\ in $\PX$ with the good rate function
\begin{equation}\label{eq-38}
    I(Q)=\inf\{i(q);q\in\PZX: q_\X=Q\}
\end{equation}
 which is (\ref{eq-35}).

\par\medskip\noindent
Let us prove (b). We rewrite the proof of Theorem \ref{res-23}
with $L^k_n$ instead of $K^k_n.$ As in the proof of Proposition
\ref{res-22}, we replace $F(z,x)$ by $f(x)$ to obtain the
pointwise convergence of the normalized log-Laplace transforms
\begin{equation}\label{eq-07}
    \limk \Lambda_k(f)=\Lambda(f)
\end{equation}
for all $f\in\CX,$ with
\begin{eqnarray*}
  \Lambda_k(f)
  &=& \frac{1}{k} \IZ \log \langle e^{kf},P^k_z\rangle\,\mu(dz) \quad
  \textrm{and}\\
  \Lambda(f)
  &=&  \IZ \supx \{f(x)-J_z(x)\}\,\mu(dz).
\end{eqnarray*}
Note that $\Lambda_k(f)=\lambda_k(F_f)$ and
$\Lambda(f)=\lambda(F_f)$ with $F_f(z,x)=f(x),$ so that
(\ref{eq-07}) is a specialization of (\ref{eq-06}).
\\
Exactly the same arguments as in the proof of Theorem \ref{res-23}
allow us to establish that $\{L^k_n\}_{k,n\geq 1}$ obeys the LDP
in $\CX'$ with  the rate function
$\Lambda^*(Q)=\sup_{f\in\CX}\left\{\langle
f,Q\rangle-\Lambda(f)\right\},$ $Q\in\CX'.$ In particular,
(\ref{eq-08}) and (\ref{eq-31}) become
\begin{equation}\label{eq-09}
     |\Lambda_k(f)|\leq \|f\|,\quad |\Lambda(f)|\leq \|f\|,
\end{equation}
for all $f\in\CX$ and
\begin{equation}\label{eq-320}
    \Glimk \Lambda_k^*=\Lambda^*
\end{equation}
in $\CX',$ where these convex conjugates are taken with respect to
the duality $(\CX',\CX).$

Thanks to (\ref{eq-217}), (\ref{eq-38}) and the uniqueness of the
rate function (see \cite{Leo05a}), we see that
$\{\Lambda^*<+\infty\}\subset \PX.$ We conclude as in the proof of
Theorem \ref{res-23} that $\{L^k_n\}_{k,n\geq 1}$ obeys the LDP in
$\PX$ with the rate function
$\Lambda^*(Q)=\sup_{f\in\CX}\left\{\langle
f,Q\rangle-\Lambda(f)\right\},$ $Q\in\PX.$ As the rate function is
unique,
\begin{equation}\label{eq-34}
    I=\Lambda^*.
\end{equation}
Considering $-f$ instead of $f$ in $\sup_{f\in\CX}$ leads to
statement (b).
\par\medskip\noindent
 Let us prove (c). As $i$ is a good rate function, the result
follows from the identity (\ref{eq-38}).
\par\medskip\noindent
Finally, statement (d) is a direct consequence of Lemma
\ref{res-310} below.
\endproof

Let us introduce the $[0,\infty]$-valued functions $I$ and $I_k$
on $\CX'$ which are defined for all $k\geq 1$ and $Q\in\CX'$ by

\begin{eqnarray}
     I_k(Q)&=& \inf\left\{\IZ \frac{1}{k} H(q^z|P^k_z)\,\mu(dz); q\in\PZX: q_\Z=\mu,q_\X=Q\right\}\label{eq-39a}\\
    I(Q)&=& \inf\left\{\IZX  J_z(x) \,q(dzdx); q\in\PZX: q_\Z=\mu,q_\X=Q
    \right\}\label{eq-39b}
\end{eqnarray}
where we use the same notation $I(Q)$ for the function on $\UX$
and its restriction to $\PX$ (see (\ref{eq-38})) and the
convention that $\inf\emptyset=+\infty.$ In particular, the
effective domains of $I_k$ and $I$ are included in $\PX.$

\par\medskip
As a by-product of the proof of Theorem \ref{res-33}, we have the
following corollary.

\begin{cor}\label{res-36}
[Hypotheses of Theorem \ref{res-33}].
 The sequence $(I_k)_{k\geq 1}$ $\Gamma$-converges to
 $I$  in $\CX'.$
\end{cor}
\proof We have shown at (\ref{eq-34}) that $\Lambda^*=I.$ It is
also true that $\Lambda_k^*=I_k,$ as can be shown by a minor
modification of the proof of (\ref{eq-210}). One concludes with
(\ref{eq-320}).
\endproof

During the proof of Theorem \ref{res-33}, we have invoked the
following

\begin{lemma}\label{res-310}
[Hypotheses of Theorem \ref{res-33}].
 If for each $k$ the Feller system
$(P_z^k)_{z\in\Z}$ satisfies (\ref{eq-301}) with $\beta$
continuous, $I$ is given by
 \begin{equation}\label{eq-42}
I(Q)=\left\{\begin{array}{ll}
     \IX J_{\beta(x)}(x)\,Q(dx)
      & \textrm{if }Q\in\PX \textrm{ and } \beta\push Q=\mu\\
      +\infty & \textrm{otherwise} \\
    \end{array}\right.,\quad Q\in\UX.
\end{equation}
\end{lemma}

\proof Let us first show that $\dom I$ is included in
$P_\beta(\mu):=\{Q\in\CX'; Q\in\PX, \beta\push Q=\mu\},$ whenever
$\beta$ is continuous.
\\
 As a direct consequence of Proposition \ref{res-22}-c, we
obtain for all $Q\in\CX'$ that
\begin{equation}\label{eq-302}
 I_k(Q)=\left\{\begin{array}{ll}
     \frac{1}{k} H(Q|P^k)
      & \textrm{if }Q\in\PX \textrm{ and } \beta\push Q=\mu\\
      +\infty & \textrm{otherwise} \\
    \end{array}\right..
    \end{equation}
This holds with $\beta$ measurable, see Remark \ref{rem-21}.
Hence,
 $\dom I_k\subset P_\beta(\mu)$ for each $k.$ Corollary
\ref{res-36} implies that $\dom I$ is included in the closure of
$P_\beta(\mu)$ in $\CX'.$ As $\beta$ is assumed to be continuous,
$\{Q\in\CX'; \langle Q,g\circ \beta\rangle=\langle\mu,g\rangle,
\forall g\in\CZ\}$ is closed in $\CX'$ and one obtains the
inclusion $\dom I\subset \{Q\in\CX'; \langle Q,g\circ
\beta\rangle=\langle\mu,g\rangle, \forall g\in\CZ\}.$ On the other
hand, $\dom I\subset\PX.$ Therefore, we obtain the desired
inclusion $\dom I\subset \{Q\in\CX'; \langle Q,g\circ
\beta\rangle=\langle\mu,g\rangle, \forall g\in\CZ\}\cap\PX=
P_\beta(\mu).$
\\
This implies that (\ref{eq-39b}) admits the unique minimizer
$q^*(dzdx)=\mu(dz)Q(dx\mid  \beta(X)=z)$ and gives (\ref{eq-42}).
\endproof

\section{Applications to the optimal transport}\label{sec-applications}

We apply the main results of Sections \ref{sec-doubly indexed} and
\ref{sec-doubly main} to the setting of Section \ref{sec-main
results}. The space $\X=\Rdd$ is the space of the random couples
and $\Z=\Rd$ is the space of the initial positions. The empirical
random measures $\Nkn$ and $\Mkn$ are specified by (\ref{eq-13a}),
(\ref{eq-04}) and (\ref{eq-03}). In the whole present section, the
Assumptions \ref{assumptions} are supposed to hold.

The spaces $\CR$ and  $\CRR$ of all continuous bounded functions
on $\Rd$ and $\Rdd$ are equipped with their topologies of uniform
convergence and their dual spaces $\CR'$ and $\CRR'$ are equipped
with the corresponding $\ast$-weak topologies, see Section
\ref{sec-conventions}.
It is convenient to use the notation $$\xi_A(y)=\left\{%
\begin{array}{ll}
    0 & \hbox{if }y\in A \\
    +\infty & \hbox{if }y\not\in A \\
\end{array}%
\right.$$ which is called the ``convex'' indicator of the subset
$A$ ($\xi_A$ is a convex function if and only if $A$ is a convex
set). Under the Assumptions \ref{assumptions}, the assumptions of
Theorem \ref{res-33} are satisfied with
$$
J_z(x)=c(x_0,x_1)+\xi_{x_0=z},\quad x=(x_0,x_1)\in \Rdd,z\in \Rd
$$
where $c$ is given at (\ref{eq-c}).  Let $\pk$ be defined by
(\ref{eq-02}). In the present setting, the functions $I_k$ and $I$
defined at (\ref{eq-39a}) and (\ref{eq-39b}) are given for all
$\rho\in \CRR'$ by $I_k=S_k$ and $I=S$ where
\begin{eqnarray*}
  S_k(\rho) &=& \frac{1}{k}H(\rho|\pk)+\xi_{\Pi_0(\mu)}(\rho) \\
  S(\rho) &=& \int_{\Rdd}c\,d\rho+\xi_{\Pi_0(\mu)}(\rho)
\end{eqnarray*}
with $\int_{\Rdd}c\,d\rho=\int_{\Rdd} c(x_0,x_1)\,\rho(dx_0dx_1)$
and
 $$
 \Pi_0(\mu)=\{\rho\in \CRR'; \langle \rho,\varphi\circ
X_0\rangle=\langle\mu,\varphi\rangle,\forall\varphi\in \CR\}
 $$
and the convention that $H(\rho|\pk)=+\infty$ and
$\int_{\Rdd}c\,d\rho=+\infty$ for all $\rho\in
\CRR'\setminus\PRR.$ Of course, $\Pi_0(\mu)\cap \PRR$ is the set
of all probability measures on ${\Rdd}$ such that $\rho_0=\mu.$
\\
The  reason for introducing $C'$ besides $\mathcal{P},$ is that
the strong unit ball $U$ of $C'$ is $\ast$-weak compact, while
compactness in $\mathcal{P}$ requires tightness criteria. This
will considerably simplify the compactness arguments.
\\
To see that the  identity about $S$ holds true, observe that the
canonical  projection $X_0$ is continuous. In particular, we have
(\ref{eq-301}) with the continuous function $\beta=X_0,$ which by
Lemma \ref{res-310} gives (\ref{eq-42}). The identity about $S_k$
is (\ref{eq-302}) with $\beta=X_0.$
\\
We shall also use the sets
\begin{eqnarray*}
  \Pi_1(\nu) &=& \{\rho\in \CRR'; \langle \rho,\varphi\circ
X_1\rangle=\langle\nu,\varphi\rangle,\forall\varphi\in \CR\}\quad \textrm{and} \\
  \Pi(\mu,\nu) &=& \Pi_0(\mu)\cap \Pi_1(\nu).
\end{eqnarray*}
As $X_0$ and $X_1$ are continuous, $\Pi_0(\mu)$ and $\Pi_1(\nu)$
are well-defined subsets of $\CRR'$ since $\varphi\circ X_0$ and
$\varphi\circ X_1$ are in $\CRR.$ We use the same notation for
$\Pi(\mu,\nu)$ in $\PRR$ and $\CRR'.$ We define for all
$\nu\in\PR$ and all $k$
\begin{eqnarray*}
  T_k(\nu) &=& \inf_{\rho\in \Pi(\mu,\nu)} \frac{1}{k}H(\rho|\pk)\\
  T(\nu) &=& \inf_{\rho\in \Pi(\mu,\nu)} \int_{\Rdd}c\,d\rho
\end{eqnarray*}
and we set $T_k(\nu)=T(\nu)=+\infty$ whenever
$\nu\in\CR'\setminus\PR.$

\noindent\emph{Caution.}\ We'll denote similarly the rate
functions $S_k, S, T_k$ and $T$ on $C'$ and their restrictions to
$\mathcal{P}.$

\begin{lemma}\label{res-12}
For each $k,$
\begin{enumerate}
    \item[(a)] $\{\Mkn\}_{n\geq 1}$ obeys the \nldp\ in $\PRR$ and
    $\CRR'$ with the good rate function $kS_k$ and
    \item[(b)] $\{\Nkn\}_{n\geq 1}$ obeys the \nldp\ in $\PR$ and
    $\CR'$ with the good rate function $kT_k.$
\end{enumerate}
\end{lemma}
\proof To get (a), apply Proposition \ref{res-22}; (b) follows by
the contraction principle.
\endproof

Applying Theorem \ref{res-33}, one obtains

\begin{thm}\label{res-10}
The following assertions hold true.
\begin{itemize}
    \item[(a)] $\{\Nkn\}_{k,n\geq}$ obeys the \knldp\ in $\PR$
    with the rate function $\nu\in\PR\mapsto \mathcal{T}_c(\mu,\nu)\in [0,\infty].$
    \item[(b)] For all $\nu\in\PR,$
    \begin{equation*}
    \mathcal{T}_c(\mu,\nu)
    =\sup_{f\in\CR}\left\{\int_{\Rd}\mathcal{S}_1f(x_0)\,\mu(dx_0)-\int_{\Rd}f(x_1)\,\nu(dx_1)\right\}
\end{equation*}
with $\mathcal{S}_1f(z)=\inf_{x_1\in\Rd}\{c(z,x_1)+f(x_1)\},$
$z\in\Rd.$
\end{itemize}
\end{thm}
\begin{remark}
The statement (b) of this theorem is the \emph{Kantorovich
duality} (\cite{Vill03}, Theorem 1.3) and Theorem \ref{res-33}-(b)
is a general version of this duality result.
\end{remark}

Similarly, we have the

\begin{prop}\label{res-11}
The following assertions hold true.
\begin{itemize}
    \item[(a)] $\{\Mkn\}_{k,n\geq}$ obeys the \knldp\ in $\PRR$
    with the rate function $S.$
    \item[(b)] For all $\rho\in\PRR$ such that
$\rho_0=\mu$
\begin{equation*}
   \int_{\Rdd}c\,d\rho
    =\sup_{g\in\CRR}\left\{\int_{\Rd}\mathcal{S}_{01}g(x_0)\,\mu(dx_0)-\int_{\Rdd}g(x_0,x_1)\,\rho(dx_0dx_1)\right\}
\end{equation*}
with $\mathcal{S}_{01}g(z)=\inf_{x_1\in\Rd}\{c(z,x_1)+g(z,x_1)\},$
$z\in\Rd.$
\end{itemize}
\end{prop}

As a consequence of the preceding results, we have the

\begin{thm}\label{res-13}
The following assertions hold true
\begin{enumerate}
    \item[(a)] $\Glimk S_k=S$ in $\CRR'$ and $\PRR$
    \item[(b)] $\Glimk T_k=T$ in $\CR'$ and $\PR.$
    \item[(c)] Since $\Glimk T_k=T,$  for all $\nu\in\PR$ there exists a
sequence $(\nu_k)$ in $\PR$ such that
     $\limk\nu_k=\nu$ in $\PR$ and
     $\limk T_k(\nu_k)=T(\nu)$ in $[0,\infty].$
\end{enumerate}
\end{thm}

\proof It is proved in \cite{Leo05a} that in a Polish space $\X,$
if one has a $k$-indexed family of \nldp s with rate functions
$kI_k$ such that the doubly indexed sequence obeys the (weak)
\knldp\ with rate function $I,$ then $\Glimk I_k=I$ in $\X.$ By
Lemma \ref{res-12} and Theorem \ref{res-10}, it follows that the
announced limits hold in the Polish spaces $\PRR$ and $\PR.$ They
also hold in $\CR'$ and $\CRR'$ since the effective domains of
$S_k$ and $S$ and of $T_k$ and $T$ (considered as functions on
$C'$) are included in $\PRR$ and $\PR.$ This proves (a) and (b).
Statement (c) follows from \cite{DalMaso}, Proposition 8.1.
\endproof

Let $\{g_i;i\geq 1\}$ be a countable subset of $\CR$ such that
$d(\gamma,\nu)=\sum_{i\geq 1}2^{-i}(|\langle
g_i,\gamma-\nu\rangle|\wedge 1),$ $\gamma,\nu\in\PR$ is a metric
which is compatible with the narrow convergence topology on $\PR.$
For all $\nu\in\PR$ and all $\rho\in \CRR',$ define
$$
d(\rho_1,\nu)=\sum_{i\geq 1}2^{-i}(|\langle g_i\circ
X_1,\rho\rangle-\langle g_i,\nu\rangle|\wedge 1).
$$
Let us recall the three minimization problems

\begin{align}
    &\textrm{minimize\quad} \frac{1}{k}H(\rho|\pk)+\alpha d(\rho_1,\nu) \textrm{\quad subject
    to\quad} \rho\in\Pi_0(\mu).\tag{MK$_k^\alpha$}\label{MKka}\\
    &\textrm{minimize\quad} \int_{\Rdd}c\,d\rho+\alpha d(\rho_1,\nu) \textrm{\quad subject
    to\quad} \rho\in\Pi_0(\mu).\tag{MK$^\alpha$}\label{MKa}\\
    &\textrm{minimize\quad} \int_{\Rdd}c\,d\rho \textrm{\quad subject
    to\quad} \rho\in\Pi(\mu,\nu).\tag{MK}\label{MK}
\end{align}

\begin{thm}\label{res-14}
Assume that $\mathcal{T}_c(\mu,\nu)<\infty.$
\begin{itemize}
    \item[(a)] We have:\quad
$\lima\limk
\inf_{\rho\in\Pi_0(\mu)}\left\{\frac{1}{k}H(\rho|\pk)+\alpha
d(\rho_1,\nu)\right\}=\mathcal{T}_c(\mu,\nu).$
    \item[(b)] For each $k$ and $\alpha,$ (\ref{MKka}) admits a unique solution
    $\rho^\alpha_k$ in $\PRR.$ For each $\alpha,$ $(\rho_k^\alpha)_{k\geq
    1}$ is a relatively compact sequence in $\PRR$ and  any limit point of $(\rho_k^\alpha)_{k\geq
    1}$ is a solution of (\ref{MKa}).
    \item[(c)] For each $\alpha,$ (\ref{MKa}) admits at least a (possibly not unique) solution $\rho^\alpha.$
    The sequence $(\rho^\alpha)_{\alpha\geq 1}$ is relatively compact in
    $\PRR$ and  any limit point of $(\rho^\alpha)_{\alpha\geq
    1}$ is a solution of (\ref{MK}).
\end{itemize}
\end{thm}

 \proof
 We introduce  functions on $\CRR'$ corresponding to (\ref{MKka}), (\ref{MKa}) and (\ref{MK}).  They are defined
for all $\rho\in \CRR'$ and each $k, \alpha\geq 1$ by
\begin{eqnarray*}
  G_k^\alpha(\rho) &=& S_k(\rho)+\alpha d(\rho_1,\nu)=\frac{1}{k}H(\rho|\pk)+\xi_{\Pi_0(\mu)}(\rho)+\alpha d(\rho_1,\nu)\\
  G^\alpha(\rho) &=& S(\rho)+\alpha d(\rho_1,\nu)=\int_{\Rdd}c\,d\rho+\xi_{\Pi_0(\mu)}(\rho)+\alpha
  d(\rho_1,\nu)\\
   G(\rho)&=&\int_{\Rdd}c\,d\rho+\xi_{\Pi(\mu,\nu)}(\rho).
\end{eqnarray*}
The domains of $S_k$ and $S$ are included in the strong unit ball
$U_{\mathbb{R}^{2d}}$ of $\CRR'.$ Therefore, the domains of
$G^\alpha_k, G_k$ and $G$ are also in $U_{\mathbb{R}^{2d}}$  which
is $\sigma(\CRR',\CRR)$-compact.
\\
We know that $S_k, S$ are \lsc, $d(\rho_1,\nu)$ is continuous and
bounded below and $\Pi_1(\nu)$ is closed. Therefore,  $G^\alpha_k,
G^\alpha$ and $G$ are inf-compact.
\\
As the relative entropy is stricly convex, $G_k^\alpha$ is also
strictly convex: it admits a unique minimizer $\rho_k^\alpha.$

As a function of $\rho,$ $d(\rho_1,\nu)$ is a finite continuous
function on $\CRR'.$ Together with the convergence $\Glimk S_k=S,$
this implies (see \cite{DalMaso}, Proposition 6.21) that for all
$\alpha,$
\begin{equation*}
    \Glimk G^\alpha_k=G^\alpha\quad \textrm{in }\CRR'.
\end{equation*}
 Observe that $\lima\alpha
d(\rho_1,\nu)=\xi_{\Pi_1(\nu)}(\rho)$ for all $\rho\in\PRR.$ As
this limit is increasing, by \cite{DalMaso}, Proposition 5.4 we
have
\begin{equation*}
    \Glima G^\alpha=G \quad \textrm{in } \CRR'.
\end{equation*}

Together with the relative compactness of the domains, these
$\Gamma$-convergence results entail the whole theorem (see
\cite{DalMaso}, Theorem 7.8 and Corollary 7.20).
\endproof

\section{$\Gamma$-convergence of convex functions on a weakly compact
space}\label{sec-appendix}

This section is dedicated to the proof of Corollary \ref{app-A2}
which is an important tool for the proof of Theorem \ref{res-23}.

A typical result about the $\Gamma$-convergence of a sequence of
convex functions $(f_n)$ is: If the sequence of the convex
conjugates $(f_n^*)$ converges in some sense, then $(f_n)$
$\Gamma$-converges. Known results of this type are usually stated
in separable reflexive Banach spaces. For instance Corollary 3.13
of H. Attouch's monograph \cite{Attouch84} is

\begin{thm}\label{app-Attouch}
Let $X$ be a separable reflexive Banach space and $(f_n)$ a
sequence of closed convex functions from $X$ into
$(-\infty,+\infty]$ satisfying the equicoerciveness assumption: $
f_n(x)\geq\alpha(\|x\|)$ for all $x\in X$ and $n\geq 1$ with
$\lim_{r\rightarrow +\infty} \alpha(r)/r=+\infty. $ Then, the
following statements are equivalent
\begin{enumerate}
    \item   $f=\mathrm{seq }X_w\textrm{-}\Gamma\textrm{-}\lim_{n\rightarrow\infty} f_n$
    \item $f^*=X_s^*\textrm{-}\Gamma\textrm{-}\lim_{n\rightarrow
    \infty}f_n^*$
    \item $\forall y\in X^*,$ $f^*(y)=\lim_{n\rightarrow\infty}f_n^*(y)$
\end{enumerate}
where $X^*$ is the dual space of $X,$ $\mathrm{seq }X_w$ refers to
the weak sequential convergence in $X$ and $X_s^*$ to the strong
convergence in $X^*.$
\end{thm}

Escaping from the reflexivity assumption is quite difficult, as
can be seen in G. Beer's monograph \cite{Beer93}.

In some applications in probability, the reflexive Banach space
setting is not as natural as it is for the usual applications of
variational convergence to PDEs. For instance when dealing with
random measures on $\X,$ the narrow topology $\sigma(\PX,
C_b(\X))$ doesn't fit the above framework since $C_b(\X)$ endowed
with the uniform topology may not be separable (unless $\X$ is
compact) and is not reflexive.

 The next result is an analogue of
Theorem \ref{app-Attouch} which agrees with applications for
random probability measures. Since we didn't find it in the
literature, we give its detailed proof.

Let $X$ and $Y$ be two vector spaces in separating duality. The
space $X$ is furnished with the weak topology $\sigma(X,Y).$

We denote $\xi_C$ the indicator function of the subset $C$ of $X$
which is defined by $\xi_C(x)=0$ if $x$ belongs to $C$ and
$\xi_C(x)=+\infty$ otherwise. Its convex conjugate is the support
function of $C:$ $\xi_C^*(y)=\sup_{x\in C}\langle x,y\rangle,$
$y\in Y.$

\begin{thm}\label{app-41}
Let $(g_n)$ be a sequence of functions on $Y$ such that
\begin{itemize}
    \item[(a)] for all $n,$ $g_n$ is a real-valued convex function on
$Y,$
    \item[(b)] $(g_n)$ converges pointwise to $g:=\limn g_n,$
    \item[(c)] $g$ is real-valued and
    \item[(d)] in restriction to any finite dimensional vector subspace $Z$ of $Y,$
    $(g_n)$ $\Gamma$-converges to $g,$ i.e. $\Glimn
    (g_n+\xi_Z)=g+\xi_Z,$ where $\xi_Z$ is the indicator
    function of $Z.$
\end{itemize}
Denote the convex conjugates on $X:$ $f_n=g^*_n$ and $f=g^*.$

If in addition,
\begin{itemize}
    \item[(e)] there exists a compact set $K\subset X$ such that
    $\dom f_n\subset K$ for all $n\geq 1$ and $\dom f\subset K$
\end{itemize}
then, $(f_n)$ $\Gamma$-converges to $f$ with respect to
$\sigma(X,Y).$
\end{thm}

\begin{remark}\label{A09}
By (\cite{DalMaso}, Proposition 5.12), under the assumption (a),
assumption (d) is implied by:
\begin{itemize}
    \item[(d')] in restriction to any finite dimensional vector subspace $Z$ of $Y,$
    $(g_n)$ is equibounded, i.e. for all $y_o\in Z,$ there exists $\delta >0$
    such that
\begin{equation*}
    \sup_{n\geq 1}\sup\{|g_n(y)|; y\in Z, |y-y_o|\leq
    \delta\}<\infty.
\end{equation*}
\end{itemize}
\end{remark}

A useful consequence of Theorem \ref{app-41} is

\begin{cor}\label{app-A2}
Let $(Y,\|\cdot\|)$ be a normed space and $X$ its topological dual
space. Let $(g_n)$ be a sequence of functions on $Y$ such that
\begin{itemize}
    \item[(a)] for all $n,$ $g_n$ is a real-valued convex function on $Y,$
    \item[(b)] $(g_n)$ converges pointwise to $g:=\limn g_n$ and
    \item[(d'')]  there exists $c>0$ such that $|g_n(y)|\leq
    c(1+\|y\|)$ for all $y\in Y$ and $n\geq 1.$
\end{itemize}
Then, $(f_n)$ $\Gamma$-converges to $f$ with respect to
$\sigma(X,Y)$ where $f_n=g^*_n$ and $f=g^*.$
\end{cor}

 \proof
Under (b), (d'') implies (c). Since the functions $g_n$ are
convex, (d'') implies that $\{g_n;n\geq 1\}$ is locally
equi-Lipschitz. Therefore (d'') implies (d') and we have (d) by
Remark \ref{A09}. Finally, (d'') implies (e) with $K=\{x\in X;
\|x\|_*\leq c\}$ where $\|x\|_*=\sup_{y, \|y\|\leq 1}\langle
x,y\rangle$ is the dual norm on $X.$ Indeed, suppose that for all
$y\in Y,$ $g(y)\leq c+c\|y\|$ and take $x\in X$ such that
$g^*(x)<+\infty.$ As for all $y,$ $\xy \leq g(y)+ g^*(x),$ we get
$|\xy|/\|y\|\leq (g^*(x)+c)/\|y\|+c.$ Letting $\|y\|$ tend to
infinity gives $\|x\|_*\leq c$ which is the announced result.
\\
The conclusion follows from Theorem \ref{app-41}.
\endproof

The proof of  Theorem \ref{app-41} is postponed after the two
preliminary Lemmas \ref{app-A1} and \ref{app-43}.

\begin{lemma}\label{app-A1}
Let $f:X\rightarrow (-\infty,+\infty]$ be a \lsc\ convex function
such that $\dom f$ is included in a compact set. Let $V$ be a
closed convex subset of $X.$

Then, if $V$ satisfies
\begin{equation}\label{A01}
    V\cap \dom f\not=\emptyset \quad \textrm{or}\quad V\cap
    \cl\dom f=\emptyset,
\end{equation}
we have
\begin{equation}\label{A02}
    \infV f(x)=-\infY (f^*(y)+\xi^*_V(-y))\in(-\infty,\infty]
\end{equation}
and if $V$ doesn't satisfy (\ref{A01}), we have
\begin{equation}\label{A03}
    \inf_{x\in W}f(x)=-\infY(f^*(y)+\xi^*_{W}(-y))=+\infty
\end{equation}
for all closed convex set $W$ such that $W\subset\inter V.$
\end{lemma}

 \proof
The proof is divided in two parts. We first consider the case
where $V\cap \dom f\not=\emptyset,$ then the case where $V\cap
\cl\dom f=\emptyset.$

$\bullet$ \textit{The case where $V\cap \dom f\not=\emptyset.$} As
$V$ is a nonempty closed convex set, its indicator function
$\xi_V$ is a closed convex function so that its biconjugate
satisfies $\xi_V^{**}=\xi_V,$ i.e. $\xi_V(x)=\supY
\{\xy-\xi_V^*(y)\}$ for all $x\in X.$ Consequently,
\begin{equation*}
    \infV f(x)=\infX\supY\{f(x)+\xy-\xi_V^*(y)\}.
\end{equation*}
One wishes to invert $\infX$ and $\supY$ by means of the following
standard inf-sup theorem (see \cite{Ekeland74} for instance). We
have $\infX\supY F(x,y)=\supY\infX F(x,y)$ provided that
$\infX\supY F(x,y)\not = \pm\infty$ and
\begin{itemize}
    \item[-] $\dom F$ is a product of convex sets,
    \item[-] $x\mapsto F(x,y)$ is convex and \lsc\ for all $y,$
    \item[-] there exists $y_o$ such that $x\mapsto F(x,y_o)$ is
    inf-compact and
    \item[-] $y\mapsto F(x,y)$ is concave for all $x.$
\end{itemize}
Our assumptions on $f$ allow us to apply this result with
$F(x,y)=f(x)+\xy-\xi_V^*(y).$ Note that
\begin{equation}\label{A04}
    \infX f(x)>-\infty
\end{equation}
  since $f$ doesn't take the value
$-\infty$ and is assumed to be \lsc\ on a compact set. Therefore,
if $\infV f(x)<+\infty,$ we have
\begin{equation*}
    \infV
    f(x)=\supY\infX\{f(x)+\xy-\xi_V^*(y)\}=-\infY\{f^*(y)+\xi_V^*(-y)\}.
\end{equation*}

 $\bullet$ \textit{The case where $V\cap \cl\dom f=\emptyset.$}
As $\cl\dom f$ is assumed to be compact, by Hahn-Banach theorem
$\cl\dom f$ and $V$ are strictly separated: there exists $y_o\in
Y$ such that $\xi_V^*(y_o)=\sup_{x\in V} \langle x,y_o\rangle <
\inf_{\cl\dom
    f}\langle x,y_o\rangle\leq \inf_{x\in\dom f}\langle
    x,y_o\rangle.$ Hence,
\begin{equation}\label{A05}
\inf_{x\in\dom f}\{\langle x,y_o\rangle-\xi_V^*(y_o)\}>0
\end{equation}
and
\begin{eqnarray*}
 -\infY (f^*(y)+\xi^*_V(-y))
 &=& \supY\infX  \{f(x)+\xy-\xi_V(y)\}\\
   &=& \supY\inf_{x\in\dom f} \{f(x)+\xy-\xi_V(y)\} \\
   &\geq& \infX f(x) +\sup_{a>0}\inf_{x\in\dom f}\{\langle x,ay_o\rangle -\xi_V^*(ay_o)\} \\
   &=& \infX f(x) +\sup_{a>0}a\inf_{x\in\dom f}\{\langle x,y_o\rangle -\xi_V^*(y_o)\} \\
   &=& +\infty
\end{eqnarray*}
where the last equality follows from (\ref{A04}) and (\ref{A05}).
This proves that (\ref{A03}) holds with $W=V.$

$\bullet$ Finally, if (\ref{A01}) isn't satisfied, taking $W$ such
that $W\subset\inter V$ insures the strict separation of $W$ and
$\cl\dom f$ as above.
\endproof

\begin{lemma}\label{app-43}
Let the $\sigma(X,Y)$-closed convex neighbourhood $V$ of the
origin be defined by
\begin{equation*}
    V=\{x\in X; \langle y_i,x\rangle\leq 1, 1\leq i\leq k\}
\end{equation*}
with $k\geq 1$ and $y_1,\dots, y_k\in Y.$ Its support function
$\xi_V^*$ is $[0,\infty]$-valued, inf-compact and its domain is
the finite dimensional convex cone spanned by $\{y_1,\dots,y_k\}.$
More precisely, its level sets are $\{\xi_V^*\leq b\}=b\, \cv
\{y_1,\dots,y_k\}$ for each $b\geq 0$ where $\cv
\{y_1,\dots,y_k\}$ is the convex hull of $\{y_1,\dots,y_k\}.$
\end{lemma}

\proof The closed convex set $V$ is the polar set of
$N=\{y_1,\dots,y_k\}:$ $V=N^\circ.$  Let $x_1\in V$ and $x_o\in
E:=\cap_{1\leq i\leq k}\mathrm{ker\,}y_i.$ Then, $\langle
y_i,x_1+x_o \rangle = \langle y_i,x_1\rangle\leq 1.$ Hence,
$x_1+x_o\in V.$ Considering the factor space $X/E,$ we now work
within a finite dimensional vector space whose algebraic dual
space is spanned by $\{y_1,\dots,y_k\}.$

We still denote by $X$ and $Y$ these finite dimensional spaces. We
are allowed to apply the finite dimension results which are proved
in the book \cite{RW98} by Rockafellar and Wets. In particular,
one knows that if $C$ is a closed convex set in $Y,$ then the
gauge function $\gamma_C(y):=\inf\{\lambda\geq 0; y\in\lambda C\},
y\in Y$ is the support function of its polar set $C^\circ=\{x\in
X; \langle x,y\rangle\leq 1, \forall y\in C\}.$ This means that
$\gamma_C=\xi_{C^\circ}^*$ (see \cite{RW98}, Example 11.19).

As $V=(N^{\circ\circ})^{\circ}$ and $N^{\circ\circ}$ is the closed
convex hull of $N,$ i.e. $N^{\circ\circ}=\cv(N):$ the convex hull
of $N,$ we get $V=\cv(N)^\circ$  and
$$
\xi_V^*=\gamma_{\cv(N)}.
$$
In particular, for all real $b,$ $\xi_V^*(y)\leq b \Leftrightarrow
\gamma_{\cv(N)}(y)\leq b \Leftrightarrow y\in b\,\cv(N).$ It
follows that the effective domain of $\xi_V^*$ is the convex cone
spanned by $y_1,\dots, y_k$ and $\xi_V^*$ is inf-compact.
\endproof

 \proof[Proof of Theorem \ref{app-41}]
 Let $\mathcal{N}(x_o)$ denote the set of
 all the neighbourhoods of $x_o\in X.$ We want to prove that
  $  \Glimn f_n(x_o):=\sup_{U\in\mathcal{N}(x_o)}\limn \inf_{x\in U}
    f_n(x)=f(x_o).$
Since $f$ is \lsc, we have $f(x_o)=\sup_{U\in\mathcal{N}(x_o)}
\inf_{x\in U} f(x),$ so that it is enough to show that for all
$U\in\mathcal{N}(x_o),$ there exists $V\in\mathcal{N}(x_o)$ such
that $V\subset U$ and
\begin{equation}\label{A06}
    \limn \infV  f_n(x)= \infV f(x).
\end{equation}
The topology $\sigma(X,Y)$ is such that $\mathcal{N}(x_o)$ admits
the sets
$$
V=\{x\in X; |\langle y_i,x-x_o\rangle|\leq 1, i\leq k\}
$$
as a base where  $(y_1,\dots,y_k), k\geq 1$ describes the
collection of all the finite families of  vectors in $Y.$
 By Lemma \ref{app-A1}, there exists such a $V\subset U$  which
 satisfies
 $$\infV f_n(x)=-\infY h_n(y)
 \textrm{ for all $n\geq 1$ and } \infV f(x)=-\infY h(y)$$ where we denote $h_n(y)=g_n(y)+\xi_V^*(-y)$
and $h(y)=g(y)+\xi_V^*(-y),$ $y\in Y.$

Let $Z$ denote the vector space spanned by $(y_1,\dots,y_k)$ and
$h^Z_n, h^Z$ the restrictions to $Z$ of $h_n$ and $h.$ For all
$y\in Y,$ we have
\begin{equation}\label{A08}
    \xi_V^*(-y)=-\langle x_o,y\rangle +
 \xi_{V-x_o}^*(-y)
\end{equation}
and by Lemma \ref{app-43}, the effective domain of $\xi_V^*$ is
$Z.$ Therefore, to prove (\ref{A06}) it remains to show that
\begin{equation}\label{A07}
    \limn\infY h_n^Z(y)=\infY h^Z(y).
\end{equation}

By assumptions (b) and (d), $(h^Z_n)$ $\Gamma$-converges and
pointwise converges to $h^Z.$ Note that this $\Gamma$-convergence
is a consequence of the lower semicontinuity of the convex
conjugate $\xi_V^*$ and Proposition 6.25 of \cite{DalMaso}.

Because of assumptions (a) and (c), $(h^Z_n)$ is also a sequence
of finite convex functions which  converges pointwise to the
finite function $h^Z.$ By (\cite{Roc70}, Theorem 10.8), $(h^Z_n)$
converges to $h^Z$ uniformly on any compact subset of $Z$ and
$h^Z$ is convex.

We now consider three cases for $x_o.$

\noindent \textit{The case where }$x_o\in\dom f.$\
 We already know that $(h^Z_n)$ $\Gamma$-converges to $h^Z.$ To
 prove (\ref{A07}), it remains to check that the sequence $(h^Z_n)$ is
 equicoercive (see \cite{DalMaso}, Theorem 7.8).
 \\
 For all $y\in Y,$ $g(y)-\langle x_o,y\rangle\geq
 -f(x_o)$ and (\ref{A08}) imply $h^Z(y)\geq -f(x_o)+\xi_{V-x_o}^*(-y).$
 Since, $-f(x_o)>-\infty$ and $\xi_{V-x_o}^*$ is inf-compact (Lemma \ref{app-43}), we obtain
 that $h^Z$ is inf-compact. As $(h^Z_n)$
converges to $h^Z$ uniformly on any compact subset of $Z,$ it
follows that $(h^Z_n)$ is equicoercive. This proves (\ref{A07}).

\noindent \textit{The case where }$x_o\in\cl\dom f.$\
 In this case, there exists $x_o'\in\dom f$ such that $V'=x_o'+(V-x_o)/2=\{x\in X; |\langle 2y_i,x-x_o'\rangle|\leq 1, i\leq
k\}\in\mathcal{N}(x_o')$ satisfies $x_o\in V'\subset V\subset U.$
One deduces from the previous case,  that (\ref{A07})  holds true
with $V'$ instead of $V.$

\noindent \textit{The case where }$x_o\not\in\cl\dom f.$\
 As $(h^Z_n)$ $\Gamma$-converges to $h^Z,$ by (\cite{Beer93},
 Proposition 1.3.5) we have $\limsup_{n\rightarrow \infty}\infY h_n^Z(y)\leq
 \infY h^Z(y).$ As $x_o\not\in\cl\dom f,$ for any small enough $V\in
\mathcal{N}(x_o),$
 $\infY h^Z(y)=-\inf_{x\in V}f(x)=-\infty.$ Therefore, $\limn\infY
 h_n^Z(y)=\infY
 h(y)=-\infty$ which is (\ref{A07}).

 This completes the proof of Theorem \ref{app-41}.
\endproof


\begin{thebibliography}{10}

\bibitem{Attouch84}
H.~Attouch.
\newblock {\em Variational convergence for functions and operators}.
\newblock Pitman Advanced Publishing Program. Pitman, 1984.

\bibitem{Beer93}
G.~Beer.
\newblock {\em Topologies on closed and closed convex sets}, volume 268 of {\em
  Mathematics and Its Applications}.
\newblock Kluwer Academic Publishers, 1993.

\bibitem{BET99}
C.~Boucher, R.S. Ellis, and B.~Turkington.
\newblock Spatializing random measures: doubly indexed processes and the large
  deviation principle.
\newblock {\em Ann. Probab.}, 27:297--324, 1999.

\bibitem{Braides}
A.~Braides.
\newblock {\em $\Gamma$-convergence for Beginners}.
\newblock Oxford Lecture Series in Mathematics 22. Oxford University Press,
  2002.

\bibitem{Brenier87}
Y.~Brenier.
\newblock D\'ecomposition polaire et r\'earrangement monotone des champs de
  vecteurs.
\newblock {\em C.R. Acad. Sci. Paris, S\'erie I}, 305:805--808, 1987.

\bibitem{CL95}
P.~Cattiaux and C.~L\'eonard.
\newblock Large deviations and {N}elson's processes.
\newblock {\em Forum Math.}, 7:95--115, 1995.

\bibitem{DG87}
D.~A. Dawson and J.~G\"artner.
\newblock Large deviations from the {M}c{K}ean-{V}lasov limit for weakly
  interacting diffusions.
\newblock {\em Stochastics}, 20:247--308, 1987.

\bibitem{DZ}
A.~Dembo and O.~Zeitouni.
\newblock {\em Large Deviations Techniques and Applications. {S}econd edition}.
\newblock Applications of Mathematics 38. Springer Verlag, 1998.

\bibitem{Ekeland74}
I.~Ekeland.
\newblock {\em La th\'eorie des jeux et ses applications \`a l'\'economie
  math\'ematique}.
\newblock Presses Universitaires de France, 1974.

\bibitem{FeyUstu04a}
D.~Feyel and A.~S. \"Ust\"unel.
\newblock {M}onge-{K}antorovitch measure transportation and {M}onge-{A}mp\`ere
  equation on {W}iener space.
\newblock {\em Probab. Theory Related Fields}, 128(3):347--385, 2004.

\bibitem{FeyUstu04b}
D.~Feyel and A.~S. \"Ust\"unel.
\newblock {M}onge-{K}antorovitch measure transportation, {M}onge-{A}mp\`ere
  equation and the {I}t\^o calculus.
\newblock In {\em Stochastic analysis and related topics in Kyoto}, volume~41
  of {\em Adv. Stud. Pure Math. Math. Soc. Japan}, pages 49--74, Tokyo, 2004.

\bibitem{Kanto42}
L.~V. Kantorovich.
\newblock On the translocation of masses.
\newblock {\em C. R. (Dokl.) Acad. Sci. URSS}, 37:199--201, 1942.

\bibitem{Kanto48}
L.~V. Kantorovich.
\newblock On a problem of {M}onge (in {R}ussian).
\newblock {\em Uspekhi Mat. Nauk.}, 3:225--226, 1948.

\bibitem{Leo05a}
C.~L\'eonard.
\newblock Large deviations of doubly indexed systems.
\newblock Preprint, 2005.

\bibitem{DalMaso}
G.~Dal Maso.
\newblock {\em An Introduction to $\Gamma$-Convergence}.
\newblock Progress in Nonlinear Differential Equations and Their Applications
  8. Birkh\"auser, 1993.

\bibitem{Mikami04}
T.~Mikami.
\newblock Monge's problem with a quadratic cost by the zero-noise limit of
  $h$-path processes.
\newblock {\em Probab. Theory Relat. Fields}, 129:245--260, 2004.

\bibitem{Monge}
G.~Monge.
\newblock M\'emoire sur la th\'eorie des d\'eblais et des remblais.
\newblock In {\em Histoire de l'Acad\'emie Royale des Sciences de Paris}, pages
  666--704. 1781.

\bibitem{Neveu}
J.~Neveu.
\newblock {\em Bases math\'ematiques du calcul des probabilit\'es}.
\newblock Masson, Paris, 1970.

\bibitem{RacRus}
S.~Rachev and L.~R\"uschendorf.
\newblock {\em Mass Transportation Problems. Vol I : Theory, Vol. II :
  Applications}.
\newblock Probability and its applications. Springer Verlag, New York, 1998.

\bibitem{Roc71}
R.~T. Rockafellar.
\newblock Convex integral functionals and duality.
\newblock In E.~Zarantonello, editor, {\em Contributions to nonlinear
  functional analysis}, pages 215--235. Academic Press, New-York, 1971.

\bibitem{Roc70}
R.T. Rockafellar.
\newblock {\em Convex Analysis}.
\newblock Princeton landmarks in mathematics. Princeton University Press,
  Princeton, N.J., 1997.
\newblock First published in the Princeton Mathematical Series in 1970.

\bibitem{RW98}
R.T. Rockafellar and R.~Wets.
\newblock {\em Variational Analysis}, volume 317 of {\em Grundlehren der
  Mathematischen Wissenschaften}.
\newblock Springer, 1998.

\bibitem{Rue78}
D.~Ruelle.
\newblock {\em Thermodynamic Formalism}.
\newblock Addison Wesley, Reading, MA, 1978.

\bibitem{Vill03}
C.~Villani.
\newblock {\em Topics in Optimal Transportation}.
\newblock Graduate Studies in Mathematics 58. American Mathematical Society,
  Providence RI, 2003.

\bibitem{Vill05}
C.~Villani.
\newblock Saint-{F}lour {L}ecture {N}otes. {O}ptimal transport, old and new.
\newblock \\ Available online via {\tt
  http://www.umpa.ens-lyon.fr/\~{}cvillani/}, 2005.

\end{thebibliography}
\end{document}